\documentclass[12pt]{amsart}

\usepackage{amsmath,amssymb,amsthm}

 \usepackage{setspace}

\theoremstyle{plain}

\theoremstyle{plain}
\theoremstyle{plain}
\newtheorem{theo}{Theorem}[section]
\newtheorem{D}[theo]{Definition}
\newtheorem{lem}[theo]{Lemma}
\newtheorem{corol}{Corollary}

\newtheorem{prop}[theo]{Proposition}

\newcommand{\Z}{\mathbb{Z}}  
  
\newcommand{\N}{\mathbb{N}}

\newcommand{\R}{\mathbb{R}}

\usepackage{amsfonts}
 
\usepackage{fancyhdr}
\usepackage{amsfonts,graphicx}
 \newtheorem{rk}{Remark }

\begin{document}

\title[Global existence for the MHD system in critical spaces]
{Global existence for the MHD system in critical spaces}

\author[H. Abidi]{Hammadi Abidi}
\address{IRMAR, Universit\'e de Rennes 1\\ Campus de
Beaulieu\\ 35~042 Rennes cedex\\ France}
\email{hamadi.abidi@univ-rennes1.fr}
\author[M. Paicu]{Marius Paicu}
\address{Laboratoire de MathŽmatique\\ 
Universit\'e Paris Sud \\
B\^atiment 425 \\
91~405 ORSAY  \\ FRANCE}
\email{marius.paicu@math.u-psud.fr}

\begin{abstract}
\noindent In this article, we show that the
magneto-hydrodynamic system (MHD) in $\R^N$ with variable density,
variable viscosity and variable conductivity has a local weak solution 
 in the Besov space $\dot
B^{\frac{N}{p_1}}_{p_1,1}(\R^N)\times\dot 
B^{\frac{N}{p_2}-1}_{p_2,1}(\R^N)
\times\dot B^{\frac{N}{p_2}-1}_{p_2,1}(\R^N)$
for all $1<p_2<+\infty$ and some
$1<p_1\leq\frac{2N}{3}$ if the initial density approaches a
positive constant. Moreover, this solution is unique if we impose the
restrictive condition $1<p_2\leq2N$. We prove also that the constructed 
solution exist globally in time if the initial data are small enough. In particular, this 
allows us to work in the frame of Besov space with negative regularity indices and this fact is
particularly important when the initial data are strong oscillating.

\end{abstract}

\maketitle
\section{ Introduction.}
In this paper we study existence and uniqueness of
solutions for the magneto-hydrodynamic system with variable viscosity
and variable density, which  describes the coupling between the
inhomogeneous Navier-Stokes system and the Maxwell equation:
$$
{\rm(MHD)}\left\{
\begin{array}{rl}
&\hspace{-0,5cm}\partial_t\rho+{\mathop{\rm div}}(\rho u)=0
\smallskip\\
&\hspace{-0,5cm}\partial_t (\rho u) 
+{\mathop{\rm div}}(\rho u\otimes u)
-\,2{\mathop{\rm div}} \big(\mu(\rho){\mathcal{M}}\big)
+\nabla\big(\Pi +\frac{B^2}{2}\big)=
\rho f+{\mathop{\rm div}}(B\otimes B)
\medskip\\
&\hspace{-0,5cm}\partial_t B-{\mathop{\rm div}}\big(\frac{\nabla B}{\sigma(\rho)}\big) 
=
B\cdot\nabla u-u\cdot\nabla B
\medskip\\
&\hspace{-0,5cm}{\mathop{\rm div}}\, u={\mathop{\rm div}}\, B=0
\medskip\\
&\hspace{-0,5cm}(\rho,u,B)_{|t=0}=(\rho_0,u_0,B_0),
\end{array}
\right.
$$
where ${\mathcal{M}}=\frac{1}{2}(\nabla u+^t\nabla u)$ is the symmetrical 
part of the gradient,
the external force $f$ is given,
$\mu(\cdot)>0$ is the viscosity of the fluid,
$\sigma(\cdot)>0$ is the conductivity  and $\Pi(t,x)$ is the pressure
in the fluid. Moreover, we  suppose that $\sigma$ and
$\mu$ are $C^\infty$ functions and that
\begin{equation}\label{viscosite-conductivite}
0<\underline\sigma\leq \frac{1}{\sigma}\leq \bar{\sigma}<\infty\quad
\hbox{and}\quad 0<\underline{\mu}\leq\mu.
\end{equation}

The homogeneous case ($\rho=\text{const.}$) of the (MHD) system
was studied by G. Duvaut and J.-L. Lions \cite{DL}. They established
local existence and uniqueness of a solution in the classical
Sobolev spaces $H^s(\R^N), \,s\geq N$. They proved also  global
existence of the solution for small initial data.

The inhomogeneous case has been studied by many authors. 
Let us mention J.-F. Gerbeau and C. Le Bris \cite{GB} and also 
B. Desjardins and C. Le Bris \cite{DB} who studied global 
existence of weak solutions of finite energy  in $\R^3$ and in 
the torus $\mathcal T^3$. On the
other hand,  local existence of strong solutions was recently
considered by H. Abidi and T. Hmidi \cite{AH}. They  proved also global existence 
of strong solutions when the initial data are 
small in some Sobolev spaces.

The principal aim of this paper is to study the strong solutions in
some Sobolev-Besov critical spaces of negative regularity index.
Working with initial data in Besov spaces of negative regularity
allows us to choose the initial velocity and the initial magnetic
field to be very irregular (even discontinuous) functions. On the other
hand, working in spaces of negative regularity allows us to prove
that the $(MHD)$ system is globally well-posed for strongly
oscillating initial data.

In the following, we suppose that the initial density verifies
$\displaystyle\inf_{x}\rho_{0}(x)> 0$ and thus, by the maximum
principle for the transport equation, we have
$\displaystyle\inf_{x}\rho(t,x)>0$. We also suppose that the density
of the fluid is a small perturbation of a constant density which we
choose to be equal to 1. This implies that we can use the transform
$a=\frac{1}{\rho}-1$ which allows us to work with the following
system:
$$
{\rm(\widetilde{MHD})}\left\{
\begin{array}{rl}
&\hspace{-0,2cm}\partial_t a+u\cdot\nabla a=0
\smallskip\\
&\hspace{-0,2cm}\partial_t  u +u\cdot\nabla u
+(1+a)\Big\{\nabla\Pi+\nabla\big(\frac{B^2}{2}\big)
-2\,{\mathop{\rm div}}
\big(\widetilde\mu(a){\mathcal{M}}\big)\Big\}=f
\\&
\hspace{9,5cm}+(1+a)B\cdot\nabla B
\medskip\\
&\hspace{-0,2cm}\partial_t B-{\mathop{\rm div}}
\big(\widetilde\sigma(a)\nabla B\big)
=B\cdot\nabla u-u\cdot\nabla B
\medskip\\
&\hspace{-0,2cm}{\mathop{\rm div}}\, u={\mathop{\rm div}}\, B=0
\medskip\\
&\hspace{-0,2cm}(a,u,B)_{|t=0}=(a_0,u_0,B_0),
\end{array}
\right.
$$
where $\widetilde \mu(a)=\mu(\frac{1}{1+a})$ and
$\widetilde\sigma(a)=\frac{1}{\sigma(\frac{1}{1+a})}$ are regular
functions.

Let us recall the theorem proved by H. Abidi et T. Hmidi in their
recent paper \cite{AH}. We denote by ${\mathcal{P}}$ the Leray projector
on the divergence free vector fields and by 
${\mathcal{Q}}=I-{\mathcal{P}}$
the projector on the gradient type vector fields. The Besov spaces
are defined in the next section.
\begin{theo}\label{MHD}\cite{AH}
Let $1<p<6.$ There exists a constant $c$ depending on $p$ and on the
functions $\mu$ and $\sigma$ such that, for $u_0,\,B_0\in\dot
B^{\frac{3}{p}-1}_{p\,1}(\R^3)$ with ${\mathop{\rm div}}\,
u_0={\mathop{\rm div}}\, B_0=0,$ $f\in L^1(\R_+;\,\dot
B^{\frac{3}{p}-1}_{p\,1}(\R^3))$ with ${\mathcal{Q}}f$ belonging to
$L^2_{loc}(\R_+;\,\dot B^{\frac{3}{p}-2}_{p\,1}(\R^3))$ and
$a_0\in\dot B^{\frac{3}{p}}_{p\,1}(\R^3)$ where
$$
\Vert a_0\Vert_{\dot B^{\frac{3}{p}}_{p\,1}}\leq c,
$$
then, there exists $T\in(0,+\infty]$ ``such that" the system
$(\rm\widetilde{MHD})$ has a solution $(a,u,B,\nabla\Pi)$
 $$
 a\in C_b\big([0,T);\,\dot B^{\frac{3}{p}}_{p\,1}\big)
\cap \widetilde L^\infty([0,T);\dot B^{\frac{3}{p}}_{p\,1})
;\,u,B\in C_b([0,T);\,\dot B^{\frac{3}{p}-1}_{p\,1})
\cap L^1(0,T;\,\dot B^{\frac{3}{p}+1}_{p\,1}).
$$ 
Moreover, there is a sufficiently
small constant $c_{1}>0$ such that, if
$$
\Vert u_0\Vert_{\dot B^{\frac{3}{p}-1}_{p\,1}}
+\Vert B_0\Vert_{\dot B^{\frac{3}{p}-1}_{p\,1}}
+\Vert f\Vert_{L^1(\R_+;\,\dot B^{\frac{3}{p}-1}_{p\,1})}
\leq
c_{1}\inf({\mu^1},\sigma^1),\quad
\hbox{with $\mu^1=\mu(1)$ et $\sigma^1=\widetilde\sigma(1),$}
$$
then $T=+\infty.$ If $1<p\leq 3,$ then
this solution is unique.
\end{theo}
This result can be easily generalized to the case of fluid evolving in 
the whole space $\R^N$.
However, the result does not provide uniqueness
for $N<p\leq 2N$, which would allow one to conclude that the system
$(\rm\widetilde{MHD})$ is globally well-posed for strongly
oscillating initial data . Addressing the issue of uniqueness is the 
principal motivation of our work.

\smallskip

\noindent In order to have a more clear idea of uniqueness, let us note that the
system $(\rm\widetilde{MHD})$ can be written as a coupled system
of a transport equation for the density and a Navier-Stokes type
equation for the couple $(u,B)$. Let us note also that the stabilizing
effect of strongly oscillating initial data  is well
known for the classical homogeneous Navier-Stokes equation.
Indeed, for the  Navier-Stokes system in the  homogeneous case
($\rho,B=\text{const.}$), i.e,
$$
{\rm({NS}_{\mu})}\left\{
\begin{array}{rl}
&\hspace{-0,5cm}\partial_t u +u\cdot\nabla u-\Delta u+\nabla\Pi=0
\\
&\hspace{-0,5cm}{\mathop{\rm div}}\, u=0
\\
&\hspace{-0,5cm}u_{|t=0}=u_0,
\end{array}
\right.
$$
it is classical to obtain global
existence and uniqueness of  solutions for small initial data in
the Besov space $\dot B^{-1+\frac{N}{p}}_{p,1}(\R^N)$ for all
$1<p<\infty$ (see \cite{CMP}). The Cannone-Meyer-Planchon result
generalizes the classical theorem by Fujita-Kato \cite{FK}, which
gives  existence and uniqueness of  solutions in the framework
of classical Sobolev spaces $\dot H^{\frac{N}{2}-1}(\R^N)$, to Besov
spaces of negative regularity index. The interest in such a result
comes from the fact that  initial data which are large in 
${\dot H}^{\frac N2-1}(\R^N)$ become small in the presence of oscillations in the
norm of the space $\dot B^{-1+\frac{N}{p}}_{p,1}$ when
$N<p<+\infty$. In particular, we find that the very fast
oscillations of the initial data stabilize the Navier-Stokes system
in the sense that the solution exists globally in time.
\begin{theo}(Cannone-Meyer-Planchon\;\cite{CMP})
Let $1<p<+\infty$ and let $u_0\in\dot B^{\frac{N}{p}-1}_{p,1}(\R^N)$
be a divergence free vector field. There then exists a time $T>0$
such that system ${\rm({NS}_{\mu})}$ has a unique solution.
$$
u\in C_b([0,T);\,\dot B^{\frac{N}{p}-1}_{p,1})\cap L^1(0,T;\,
\dot B^{\frac{N}{p}+1}_{p,1}).
$$
Moreover, there is a constant $c>0$ small enough such that if
$$
\Vert u_0\Vert_{\dot B^{\frac{N}{p}-1}_{p,1}}\leq c\mu,
$$
then  $T=\infty.$
\end{theo}

\smallskip
In this article we will show the existence and uniqueness of 
global solution for system $(\rm\widetilde{MHD})$ for strongly oscillating initial data. 
For that it will be necessary to work in spaces with negative index of regularity. 
Let us note that the result of \cite{AH} does not make
it possible to construct a unique global solution for the data in
spaces of negative index, since one has  uniqueness of the solution
only in the case when $1<p\leq N.$ Also let us note that one has
existence of a global weak solution when $N<p<2N$ for small data. 
In this paper, we prove in fact that the $(\rm\widetilde{MHD})$ 
system is globally well-posed for oscillating initial data, when
$\frac{1}{\rho_0}-1\in\dot B^{\frac{N}{p_1}}_{p_1,1}$ and
$u_0,\,B_0\in \dot B^{\frac{N}{p_2}-1}_{p_2,1}$ with 
$p_1\leq p_2$ and $\frac{1}{p_1}\leq \frac{1}{p_2}+\frac{1}{N}$ 
and $\frac{1}{p_1}+\frac{1}{p_2}\geq\frac{2}{N}.$ 
Note in particular that we
obtain the H. Abidi and T. Hmidi results as a particular case of our
theorem by taking $p_1=p_2.$ The improvement obtained in our result
is due directly to the fact that we work with the density, velocity field, and
magnetic field in the spaces of Besov built on different spaces of
Lebesgue. The method of the proof is based on the regularizing
effect for the heat equation (for more precise details, see
\cite{CH}). To be more precise, we point out a result of harmonic
analysis due to R. Danchin \cite{DAN}, which is an inequality of the
type of the Poincar\'e-type inequality for  functions localized in
frequencies. That enables us to gain two derivatives of the solution
from the heat equation starting from the Laplacian, and thus, for
initial data in $\dot B^{-1+\frac{N}{p}}_{p,1}(\R^N)$ we find that
the solution belongs to the space $L^1([0,T];\, \dot
B^{1+\frac{N}{p}}_{p,1})$ which is a subspace of $L^1(Lip\,(\R^N))$.
This is the principal reason for why one cannot work with the
initial data $u_0\in\dot B^{-1+\frac{N}{p}}_{p, r}$ for $r>1.$

We prove an existence result in critical Besov spaces (for the
definition see the next section). Our principal result is as
follows:
\begin{theo}\label{Res-Pri}
Let $1<p_1\leq p_2<+\infty$ be such that
$\frac{1}{p_1}\leq\frac{1}{p_2}+\frac{1}{N}$ 
and 
$\frac{1}{N}<\frac{1}{p_1}+\frac{1}{p_2}.$ 
There exists a positive constant
$c$ depending on $p$ and on functions $\mu,$ $\sigma$ such that, for
$u_0,\,B_0\in\dot B^{\frac{N}{p_2}-1}_{p_2,1}(\R^N)$ with
${\mathop{\rm div}}\, u_0={\mathop{\rm div}}\, B_0=0,$ $f\in
L^1_{loc}(\R_+;\,\dot B^{\frac{N}{p_2}-1}_{p_2,1}(\R^N))$ with
${\mathcal{Q}}f\in L^2_{loc}(\R_+;\,\dot
B^{\frac{N}{p_2}-2}_{p_2,1}(\R^N))$ and $a_0\in\dot
B^{\frac{N}{p_1}}_{p_1,1}(\R^N)$ where
$$
\Vert a_0\Vert_{\dot B^{\frac{N}{p_1}}_{p_1,1}}\leq c,
$$
then there exists $T(u_0,B_0,f)>0$ such that the system $(\widetilde{MHD})$
has a solution $(a,u,B,\nabla\Pi)$ with
 $$
 a\in C_b\big([0,T);\,\dot B^{\frac{N}{p_1}}_{p_1,1}\big)
\cap \widetilde L^\infty([0,T);\dot B^{\frac{N}{p_1}}_{p_1,1})
;\,u,B\in C_b([0,T);\,\dot B^{\frac{N}{p_2}-1}_{p_2,1})\cap
L^1(0,T;\, \dot B^{\frac{N}{p_2}+1}_{p_2,1})
$$
$$
\mbox{and}\quad
\nabla\Pi\in L^{\frac{2}{2-\eta}}_{T}(\dot B^{\frac{N}{p_2}-1-\eta}_{p_2,1}),\;\;\hbox{with  }\,
0\leq\eta<\inf(1,\frac{2N}{p_2})\;\mbox{and}\;
\frac{1}{N}+\frac{\eta}{N}<\frac{1}{p_1}+\frac{1}{p_2}.
$$
Moreover, there exists a positive constant $c_{1}$ such that if
$$
\Vert u_0\Vert_{\dot B^{\frac{N}{p_2}-1}_{p_2,1}}+\Vert B_0\Vert_{\dot
B^{\frac{N}{p_2}-1}_{p_2,1}}
+\Vert f\Vert_{L^1(\R_+;\,\dot B^{\frac{N}{p_2}-1}_{p_2,1})}
\leq
c_{1}\inf({\mu^1},\sigma^1),
$$
with $\mu^1=\widetilde\mu(1),$ $\sigma^1=\widetilde\sigma(1)$, then
$T=+\infty.$\\
If, in addition, we have that $1<p_2\leq 2N,$ and 
$\frac{1}{p_1}+\frac{1}{p_2}\geq
\frac{2}{N}$ then such a solution is unique.
\end{theo}

\smallskip

The proof of Theorem \ref{Res-Pri}, is carried out in two stages. 
Firstly, we show the uniqueness result that is based on a logarithmic
estimate combined Osgood lemma. Secondly, for the existence part we proceed as follows: we regularize both initial data and 
$(\rm\widetilde{MHD})$ system, for which we establish the existence of solutions. After we shaw we can bound from below the time existence. Finally we prove that the regularization solutions converge to a solution satisfying our initial problem.

\begin{rk}
In The case of variable viscosity and variable 
conducitvity, we need the  more restrictive condition $p_1\leq p_2$. This 
condition does not appear in the case where the viscosity is constant (see 
our paper \cite{Abidi-Paicu}).
\end{rk}

\begin{rk}
This theorem allows us to construct a solution
(local in time in general, respectively global in time when the initial 
data is small compared with viscosity), for $u_0,\,B_0\in\dot
B^{-1+\frac{N}{p_2}}_{p_2,1}(\R^N)$ and all $1<p_2<+\infty$. In
fact, it is enough for example to consider the density such that
$a_0=\rho_0^{-1}-1\in\dot B^{1}_{N\,1}(\R^N)$ when $N\leq
p_2<+\infty.$ In the case  when $1<p_2<N$ we take for example $p_1=p_2$
 (other choices are possible, it suffices for example to
take $p_1$ which verifies $\sup(1,\frac{Np_2}{N+p_2})<
p_1\leq p_2$).

\smallskip

\noindent On the other hand, we obtain a unique solution for all
$u_0,\,B_0\in\dot B^{-1+\frac{N}{p_2}}_{p_2,1}(\R^N)$ for all
$1~<~p~_2\leq ~ 2N$.  In order to obtain this, it suffices to
consider for example $a_0=\rho_0^{-1}-1\in\dot
B^{\frac{N}{p_1}}_{p_1,1}(\R^N)$ with $p_1=\frac{2N}{3}$ when $N\leq
p_2\leq 2N$, and and it suffices to take  $\sup(1,\frac{Np_2}{N+p_2})<p_1\leq p_2$
when $1<p_2<N.$ 
\end{rk}

\smallskip

\begin{rk} In particular, Theorem \ref{Res-Pri} implies
existence of a unique global solution for the
$(\widetilde{MHD})$ system, when the initial data $(\rho_0,u_0,B_0)$
have the particular form
$$
\begin{aligned}
a_0=\rho_0^{-1}-1\in&\mathcal S(\R^3); 
\hskip0.5cm 
u_0=\varepsilon^{-\alpha}\sin\bigg(\frac{x_3}{\varepsilon}\bigg)(-\partial_2\phi^1,\partial_1\phi^1,0);
\\&
\hskip0.5cmB_0=
\varepsilon^{-\beta}\sin\bigg(\frac{x_3}{\varepsilon}\bigg)
(-\partial_2\phi^2,\partial_1\phi^2,0)
\end{aligned}
$$
with $\alpha,\beta\in[0,1)$, $\inf\limits_{x\in\R^3}\rho_0>0$ and
$\phi^i\in \mathcal S(\R^3),$ with $a_0$ small and $\varepsilon>0$
small enough. Indeed, it is easy to verify the following assertion.
Let $\phi\in{\mathcal{S}}(\R^3)$, $k\in \R^3$, $|k|\neq 0$ and
$(\sigma,p,r)\in\R_+^*\times[1,\infty)^2.$ Then, the function
$\phi_\varepsilon(x)~=~\phi(x)\,e^{ix\cdot k/\varepsilon}$ is small
in the space $\dot B^{-\sigma}_{p,r}.$ More precisely, we have
$$
\Vert \phi_{\varepsilon}\Vert_{\dot B^{-\sigma}_{p,r}} \leq
C(\phi){\varepsilon}^\sigma,
$$
where $C(\phi)=\|\phi\|_{\dot B^{\sigma}_{p, r}}.$

\end{rk}

\section{ Preliminaries.}
\subsection{ Notation.}
Let $X$ be a Banach space and $p\in[1,\infty]$. We denote
by $L^p(0,T;\,X)$ the set of measurable functions
$f:(0,T)\rightarrow X$, such that $t\longmapsto\Vert f(t)\Vert_X$
belongs to $L^p(0,T),$ and we denote by $C([0,T);\,X)$ we denote the space of
continuous functions on $[0,T)$ with values in $X$, $C_b([0,T);\,
X):=C([0,T);\,X)\cap L^{\infty}(0,T;\,X).$ Let
$\mu^1=\mu(1),$ $\widetilde\mu(a)=\mu(\frac{1}{1+a}),$
$\widetilde\sigma(a)=\frac{1}{\sigma(\frac{1}{1+a})},$
$\sigma^1=\widetilde\sigma(1)$ and for $1\leq p\leq\infty,$ we
denote by $p'$ the conjugate exponent of $p$ given by
$\frac{1}{p}+\frac{1}{p'}=1.$
\subsection{Littlewood-Paley theory.}
In this section, we briefly recall the Littlewood-Paley
theory and we define the functional spaces in which we will work. To
this order, we use a unit dyadic (see for
example \cite{C}). Let $\mathcal C\subset \R^N$ be the annulus centered 
in $0$,
with the small radius $\frac{3}{4},$ and the big radius $\frac{8}{3}.$
There exist two positive radially symmetric functions $\chi$ and 
$\varphi$
belonging respectively to $C^\infty_0\big(B(0,\frac{4}{3})\big)$ and
to $C^\infty_0(\mathcal C)$ such that:
$$
\sum_{q\in \Z}\varphi(2^{-q}\xi)=1\quad\quad \forall\xi\neq
0\quad\hbox{et}\quad
\chi(\xi)+\sum_{q\in
\N}\varphi(2^{-q}\xi)=1\quad\quad\forall\xi\in\R^N.
$$
We define the following operators.
$$
\begin{aligned}
\Delta_q\;u =\varphi(2^{-q}D)\;u\quad\forall\; q\in\Z
\quad\hbox{et}\quad
S_q\;u=\sum_{p\leq q-1}\Delta_{p}v \quad\forall\; q\in\Z.
\end{aligned}
$$
Moreover, we have:
$$
u=\sum_{q\in\Z}\Delta_q \,u\quad\forall\,u\in {\mathcal
{S}}'(\R^N)/{\mathcal{P}}[\R^N],
$$
where ${\mathcal{P}}[\R^N]$ is the set of polynomials (see for example
\cite{PE}). Moreover, the Littlewood-Paley decomposition satisfies the
property of almost orthogonality:
\begin{equation}\label{Pres_orth}
\Delta_k\Delta_q u\equiv 0
\quad\mbox{if}\quad\vert k-q\vert\geq 2
\quad\mbox{and}\quad\Delta_k(S_{q-1}u\Delta_q u)
\equiv 0\quad\mbox{if}\quad\vert k-q\vert\geq 5.
\end{equation}
\begin{D}\label{def}
For $s\in\R,\,(p,r)\in[1,+\infty]^2$ and $u\in {\mathcal {S}}'(\R^N),$
we denote
$$
\Vert u\Vert_{\dot B^s_{p\,r}}:=
\Big(\sum_{q\in\Z}2^{rqs}\Vert\Delta_q
\,u\Vert_{L^p}^r\Big)^{\frac{1}{r}}
$$
with the usual change for the case $r=+\infty$. Then for
$s<\frac{N}{p}$ and $s\leq\frac{N}{p},\,r=1$ we define
$$
\dot B^s_{p\,r}:=
\Big\{u\in {\mathcal S}'(\R^N)\;\Big|\;\Vert u\Vert_{\dot
B^s_{p\,r}}<\infty\Big\},
$$
otherwise, we define $\dot B^s_{p\,r}$ like the adherence in
${\mathcal{S}}'$ of functions belonging to the Schwartz space, for the
norm $\Vert\cdot\Vert_{\dot B^s_{p\,r}}.$
\end{D}
Let us recall also the Bernstein inequality (see for example
\cite{C}) which allows us to obtain some embeddings of spaces.
\begin{lem}\label{bernstein2} (BERNSTEIN)\,
{\it Let $(r_1,r_2)$ be a couple of nonnegative real numbers such that
$r_1<r_2.$ Then there exists a nonnegative constant $C$ such that for
any integer $k,$ any couple $(a,b)$ such that $1\leq a\leq
b\leq\infty$ and every function $u$ in $L^{a}(\R^N),$ we have
$$
\begin{aligned}
&\quad\quad\mbox{Supp}\;{{\mathcal{F}}}u\in B(0,\lambda r_1)
\Longrightarrow\displaystyle\sup_{\vert\alpha\vert=k}\Vert\partial^{\alpha}u\Vert_{L^b}
\leq C^k{\lambda}^{k+N(\frac{1}{a}-\frac{1}{b})}\Vert u\Vert_{L^{a}},
\\&
\mbox{Supp}\;{{\mathcal{F}}}u\in C(0,\lambda r_1,\lambda r_2)
\Longrightarrow C^{-k}{\lambda}^k\Vert u\Vert_{L^{a}}
\leq\displaystyle\sup_{\vert\alpha\vert=k}\Vert\partial^{\alpha}u\Vert_{L^a}\leq
C^k{\lambda}^{k}
\Vert u\Vert_{L^a}.
\end{aligned}
$$}
\end{lem}
In  order to obtain a better description of the regularizing effect
of the transport-diffusion equation, we will use the spaces
$\widetilde{L}^{\rho}_T(\dot B^s_{p\,r})$ introduced by J.-Y. Chemin
and N. Lerner in \cite{CL}.
\begin{D}\label{chaleur+}
Let $s\leq\frac{N}{p}$ (respectively $s\in\R$),
$(r,\rho,p)\in[1,\,+\infty]^3$ and $T\in]0,\,+\infty]$. We say then
that $f\in \widetilde{L}^{\rho}_T(\dot B^s_{p\,r}),$ if
$$
\Vert f\Vert_{\widetilde{L}^{\rho}_T(\dot B^s_{p\,r})}
:=\Big(\sum_{q\in\Z}2^{qrs}
\Big(\int_0^T\Vert\Delta_q\,f(t)\Vert_{L^p}^{\rho}dt\Big)^{\frac{r}{\rho}}\Big)^{\frac{1}{r}}
<\infty.
$$
with the usual change if $r=\infty.$
\end{D}
\noindent For $\theta\in[0,1],$ we have
\begin{equation}\label{Interpo}
\Vert u\Vert_{\widetilde{L}^{\rho}_T(\dot B^s_{p\,r})} \leq \Vert
u\Vert_{\widetilde{L}^{\rho_1}_T( \dot B^{s_1}_{p\,r})}^{\theta}
\Vert u\Vert_{\widetilde{L}^{\rho_2}_T(\dot
B^{s_2}_{p\,r})}^{1-\theta}
\end{equation}
with $\frac{1}{\rho}=\frac{\theta}{\rho_1}+\frac{1-\theta}{\rho_2}$
and $s=\theta s_1+(1-\theta)s_2.$\\
Note that the Minkowski inequality implies that
$$
\Vert u\Vert_{\widetilde{L}^{\rho}_T(\dot B^s_{p\,r})}
\leq
\Vert u\Vert_{L^{\rho}_T(\dot B^s_{p\,r})}
\quad\mbox{if}\quad\rho\leq r
\quad\hbox{and}\quad
\Vert u\Vert_{L^{\rho}_T(\dot B^s_{p\,r})}
\leq
\Vert u\Vert_{\widetilde{L}^{\rho}_T(\dot B^s_{p\,r})}
\quad\mbox{if}\quad r\leq\rho.
$$
We give now the product laws in Besov spaces based on different
Lebesgue spaces. This product laws are studied in detail in the
paper \cite{Abidi-Paicu}.
\begin{prop}\label{Loi-Pro-Par}
Let $(p,p_1,p_2,r,\lambda_1,\lambda_2)\in[1,\infty]^6$ such that
$\frac{1}{p}\leq\frac{1}{p_1}+\frac{1}{p_2},$ $p_1\leq\lambda_2,$
$p_2\leq\lambda_1,$
$\frac{1}{p}\leq\frac{1}{p_1}+\frac{1}{\lambda_1}\leq 1$ et
$\frac{1}{p}\leq\frac{1}{p_2}+\frac{1}{\lambda_2}\leq1.$
Then, we have the following inequality:\\
If $s_1+s_2+N\inf(0,1-\frac{1}{p_1}-\frac{1}{p_2})>0,$
$s_1+\frac{N}{\lambda_2}<\frac{N}{p_1}$ and
$s_2+\frac{N}{\lambda_1}<\frac{N}{p_2}.$ Then
\begin{equation}\label{p_1p_2}
\|uv\|_{\dot 
B^{s_1+s_2-N(\frac{1}{p_1}+\frac{1}{p_2}-\frac{1}{p})}_{p,r}}
\lesssim \|u\|_{\dot B^{s_1}_{p_1,r}}\|v\|_{\dot B^{s_2}_{p_2,\infty}},
\end{equation}
when $s_1+\frac{N}{\lambda_2}=\frac{N}{p_1}$ (respectively
$s_2+\frac{N}{\lambda_1}=\frac{N}{p_2}$) we replace
$\|u\|_{\dot B^{s_1}_{p_1,r}}\|v\|_{\dot B^{s_2}_{p_2,\infty}}$ (respectively
$\|v\|_{\dot B^{s_2}_{p_2,\infty}}$) by
$\|u\|_{\dot B^{s_1}_{p_1,1}}\|v\|_{\dot B^{s_2}_{p_2,r}}$ (respectively
$\|v\|_{\dot B^{s_2}_{p_2,\infty}\cap L^\infty}$), if
$s_1+\frac{N}{\lambda_2}=\frac{N}{p_1}$ and
$s_2+\frac{N}{\lambda_1}=\frac{N}{p_2}$ we take $r=1.$
\bigskip

\noindent If $s_1+s_2=0,$
$s_1\in(\frac{N}{\lambda_1}-\frac{N}{p_2},\frac{N}{p_1}-\frac{N}{\lambda_2}]$
and $\frac{1}{p_1}+\frac{1}{p_2}\leq1,$ then
\begin{equation}\label{sommenulle}
\|uv\|_{\dot 
B^{-N(\frac{1}{p_1}+\frac{1}{p_2}-\frac{1}{p})}_{p,\infty}}
\lesssim \|u\|_{\dot B^{s_1}_{p_1,1}}\|v\|_{\dot B^{s_2}_{p_2,\infty}}.
\end{equation}
If $\vert s\vert<\frac{N}{p}$ for $p\geq2$ and
$-\frac{N}{p'}<s<\frac{N}{p}$ otherwise, we have
\begin{equation}\label{produit3}
\Vert uv\Vert_{\dot B^s_{p,r}} \lesssim \Vert u\Vert_{\dot B^s_{p,r}} 
\Vert
v\Vert_ {\dot B^{\frac{N}{p}}_{p,\infty}\cap L^{\infty}}.
\end{equation}
\end{prop}
\begin{rk}
In the following, $p$ will be equal to $p_1$ or to $p_2$ and
$\frac{1}{\lambda}=\frac{1}{p_1}-\frac{1}{p_2}$ if $p_1\leq p_2$,
respectively $\frac{1}{\lambda}=\frac{1}{p_2}-\frac{1}{p_1}$ if
$p_2\leq p_1$.
\end{rk}

\begin{rk}
Note that for $p_1=p_2$ we obtain the classical product laws.
On the other hand, if $s_i<\frac{N}{p_i}$, $s_1+s_2>0$ and $p_1\leq
p_2$ we obtain that $uv\in\dot B^{s_1+s_2-\frac{N}{p_1}}_{p_2,1}$,
otherwise, if $s_i<\frac{N}{p_2}$ we obtain $uv\in
\dot B^{s_1+s_2-\frac{N}{p_2}}_{p_1,1}$. The interpretation of this
facts, is that in a product law we can a smaller number of
derivatives than usual, if we measure these derivatives with a $L^p$
Lebesque space with small $p\geq 1.$
\end{rk}

\begin{rk}
The Proposition \ref{Loi-Pro-Par} also holds in $\widetilde 
L^\rho_t(\dot B^s_{p,r}).$ For example inequality (\ref{produit3}) becomes
$$
\Vert uv\Vert_{\dot B^s_{p,r}} \lesssim \Vert u\Vert_{\dot B^s_{p,r}} 
\Vert
v\Vert_ {\dot B^{\frac{N}{p}}_{p,\infty}\cap L^{\infty}}
$$
whenever $\vert s\vert<\frac{N}{p}$ for $p\geq2$ and
$-\frac{N}{p'}<s<\frac{N}{p},$
$1\leq\rho,\rho_1,\rho_2\leq\infty$ and
$1/\rho=1/{\rho_1}+1/{\rho_1}.$
\end{rk}
\section{Estimates for the  transport and Stokes equations.}
We note that the MHD system with variable density consists of a
transport equation for the density and a Stokes equation for the
velocity vector-field. We begin  by  giving the necessary
estimates for the transport and for the non-stationary Stokes
equations (for the proofs, see the paper \cite{Abidi-Paicu}):
\begin{prop}\label{eqtransport}
Let $(p_1,p_2)\in[1,+\infty]^2,$
$-1-N\inf({1\over p_2},{1\over p'_1})<s
<1+N\inf(\frac{1}{p_1},\frac{1}{p_2})$
where  $p'_1$ is the conjugate
exponent of $p_1$ (respectively
$s=1+N\inf(\frac{1}{p_1},\frac{1}{p_2})$) and $r\in[1, +\infty]$
(respectively $r=1$). Let $u$ a free-divergence vector field such that
$\nabla u\in L^1(0,T;\,\dot B^{\frac{N}{p_2}}_{p_2,r}\cap L^\infty)$
(respectively $u\in L^1(0,T;\,\dot B^{\frac{N}{p_2}+1}_{p_2,1})$). We
suppose that $\rho_0\in \dot B^s_{p_1,r},$ $f\in L^1(0,T;\,\dot
B^s_{p_1,r}).$ Let $\rho\in L^\infty(0,T;\,\dot B^s_{p_1,r})\cap
C([0,T];\,{\mathcal S}')$ be a solution of the following system
$$
\left\{
\begin{array}{rl}
&\partial_ t \rho+u\cdot\nabla \rho=f,\\
&\rho_{|t=0}=\rho_0.
\end{array}
\right.
$$
Then there exists a non-negative constant $C$ depending on  $N$ and $s$
such that
\begin{equation}\label{es.eq.tr}
\Vert\rho\Vert_{\widetilde{L}^{\infty}_T(\dot B^s_{p_1,r})} \leq
e^{CU(t)}\Big(\Vert\rho_0\Vert_{\dot B^s_{p_1,r}}+\,\int_0^t \Vert
f(\tau)\Vert_{\dot B^s_{p_1,r}}d\tau\Big),
\end{equation}
where $U(t)=\,\int_0^t\Vert\nabla u(\tau)\Vert_{\dot
B^{\frac{N}{p_2}}_{p_2,r}\cap L^{\infty}}d\tau.$ (respectively
$U(t)=\,\int_0^t\Vert u(\tau)\Vert_{\dot
B^{\frac{N}{p_2}+1}_{p_2,1}}d\tau$).
\end{prop}

\begin{prop}\label{eq.stokes}
Let $p\in]1,\infty[$ and
$-1-N\inf({1\over p},{1\over p'})<s<{N\over p},$
where  $p$ is the conjugate exponent of $p.$
Let $u_0$ be a divergence free
vector field with the components in  $\dot B^s_{p,r}$ and $g$ a
vector field with the components in 
$\widetilde{L}^1_T(\dot B^s_{p,r}).$ Let $u$ and $v$ be two divergence free vector fields
 such that  $\nabla v$
 has the coefficients in 
 $L^1(0,T;\,\dot B^{\frac{N}{p}}_{p,r}\cap L^\infty)$ 
 (respectively $L^1_T(\dot B^{\frac{N}{p}}_{p,1})$) and 
$u\in C([0,T;\,\dot B^s_{p,r})\cap 
\widetilde{L}^1_T(\dot B^{s+2}_{p,r}).$
Let $u$ be a solution of the non stationary Stokes system
$$
{\rm(L)}\quad\quad \left\{\begin{array}{rl}
&\partial_t u +v\cdot\nabla u-\nu\Delta u+\nabla\Pi = g \\
&{\mathop{\rm div}}\, u=0\\
&u_{|t=0}=u_0.
\end{array}
\right.
$$
Then there exists $C>0$ depending on $N$ and $s$ such that $u$ verifies
 the following estimate
 \begin{equation}\label{es.eq.st}
\Vert u\Vert_{\widetilde{L}^{\infty}_T(\dot B^s_{p,r})} 
+\nu \Vert u\Vert_{\widetilde{L}^1_T(\dot B^{s+2}_{p,r})}
+\Vert\nabla\Pi\Vert_{\widetilde{L}^1_T(\dot B^s_{p,r})}
 \leq
e^{C\Vert\nabla v\Vert_{L^1_T(\dot B^{\frac{N}{p}}_{p,r}\cap
L^\infty)}} \Big\{\Vert u_0\Vert_{\dot B^s_{p,r}}
+C\Vert g\Vert_{\widetilde{L}^1_T(\dot B^s_{p,r})}\Big\}.
\end{equation}
Moreover, if $2\leq p$ and $s=-1-{N\over p},$ then we have the
following estimate:
\begin{equation}\label{es.st.li}
\Vert u\Vert_{L^\infty_T(\dot B^{s}_{p,\infty})} +\nu\Vert
u\Vert_{\widetilde{L}^1_T(\dot B^{2+s}_{p,\infty})}
+\Vert\nabla\Pi\Vert_{\widetilde{L}^1_T(\dot B^{s}_{p,\infty})} \leq
e^{C\Vert\nabla v\Vert_{L^1_T(\dot B^{\frac{N}{p}}_{p,1})}} \Big\{\Vert
u_0\Vert_{\dot B^{s}_{p,\infty}} +C\Vert
g\Vert_{\widetilde{L}^1_T(\dot B^{s}_{p,\infty})}\Big\}.
\end{equation}
\end{prop}

\

Let us recall the Osgood Lemma (see \cite{fleet}), which allows us to
infer uniqueness of the solution in the critical case (see
the uniqueness section).

\begin{lem}(Osgood)
\label{osgood}

{\it Let $\rho\geq 0$ be a measurable function, $\gamma$ be a locally
integrable function and $\mu$ be a positive, continuous
 and non decreasing function which verifies the following condition
$$\int_0^1\frac{dr}{\mu (r)}=+\infty.$$

\noindent Let also $a$ be a positive real number and let $\rho$ satisfy
 the inequality

$$\rho (t)\leq a +\int _0^t\gamma (s)\mu(\rho (s))ds.$$

\noindent Then if $a$ is equal to zero, the function $\rho$  vanishes.

\noindent If  $a$ is not zero, then we have
$$
-\mathcal M (\rho (t))+\mathcal M (a)\leq \int_0^t \gamma (s)ds,
\hskip0.5cm\text{with}\hskip0.5cm \mathcal M
(x)=\int\limits_x^1\frac{dr}{\mu (r)}\cdot
$$}
\end{lem}

Finally, we recall the following result of logarithmic interpolation
(see \cite{D} Proposition 2.8).
\begin{lem}\label{DANCHIN}
Let $(p,\lambda)\in[1,+\infty]^2,$ $s\in\R,$ $t\in\R_+,$ 
$\varepsilon\in(0,1]$
and $u\in\widetilde L^\lambda_t(\dot B^{s-\varepsilon}_{p,\infty})
\cap \widetilde L^\lambda_t(\dot B^{s}_{p,1})
\cap \widetilde L^\lambda_t(\dot B^{s+\varepsilon}_{p,\infty}).$
Then
$$
\|u\|_{\widetilde L^\lambda_t(\dot B^{s}_{p,1})}
\lesssim
{\|u\|_{\widetilde L^\lambda_t(\dot B^{s}_{p,\infty})}\over\varepsilon}
\log
\Big(e+{\|u\|_{\widetilde L^\lambda_t(\dot B^{s-\varepsilon}_{p,\infty})}
+\|u\|_{\widetilde L^\lambda_t(\dot B^{s+\varepsilon}_{p,\infty})}
\over\|u\|_{\widetilde L^\lambda_t(\dot B^{s}_{p,\infty})}}\Big).
$$
\end{lem}
\section{ Proof of the Theorem \ref{Res-Pri}.}
We will proceed in two steps. First we prove the uniqueness of the
solution which is principally based on a logarithmic estimate and on
the Osgood Lemma which is useful in the case of logarithmic
estimates. The second part is devoted  to the proof of existence of
the solution.
\subsection{Uniqueness.}
Let $1\leq p_2\leq 2N$ and $1<p_1\leq p_2$ be such that
$\frac{1}{p_1}+\frac{1}{p_2}\geq\frac{2}{N}$ and
$\frac{1}{p_1}\leq\frac{1}{p_2}+\frac{1}{N}$. We denote by
$(a^i,u^i,\nabla\Pi^i)$ for $1\leq i\leq 2$ two solutions of the
$(\rm\widetilde{MHD})$ system. We define
$$
\begin{aligned}
&({\mathcal{M}}^i,\delta{\mathcal{M}}):=(\frac{1}{2}(\nabla u^i+^t\nabla 
u^i),{\mathcal{M}}^2-{\mathcal{M}}^1)\\&
\mbox{and}\;\,
(\delta a,\delta u,\nabla\delta\Pi,\delta B):=(a^2-a^1,u^2-u^1,\nabla\Pi^2-\nabla\Pi^1,B^2-B^1).
\end{aligned}
$$
We can easily check that
$$
\left\{
\begin{array}{rl}
&\partial_t \delta a +u^2\cdot\nabla\delta a= -\delta u\cdot\nabla a^1
\medskip\\
&\partial_t\delta u+u^2\cdot\nabla\delta u-\mu^1\Delta\delta
u+\nabla\delta\Pi
=H(a^{i},u^{i},\nabla\Pi^{i},B^i)
\medskip\\
&\partial_t\delta B+u^2\cdot\nabla\delta B-\sigma^1\Delta\delta
B=G(a^i,u^i,B^i)
\medskip\\
&{\mathop{\rm div}}\, \delta u={\mathop{\rm div}}\, \delta B=0,
\end{array}
\right.
$$
where
$$
\begin{aligned}
&H(a^{i},u^{i},\nabla\Pi^{i},B^i)
=
-\delta u\cdot\nabla u^1+a^1(\mu^1\Delta\delta u-\nabla\delta\Pi)
+\delta a(\mu^1\Delta u^2-\nabla\Pi^2)
\\&
+2{\mathop{\rm 
div}}\Big[\big(\widetilde\mu(a^2)-\mu^1\big)\delta{\mathcal{M}}\Big]
+2\delta a\,{\mathop{\rm 
div}}\Big[\big(\widetilde\mu(a^2)-\mu^1\big){\mathcal{M}}^2\Big]
+2a^1{\mathop{\rm 
div}}\Big[\big(\widetilde\mu(a^1)-\mu^1\big)\delta{\mathcal{M}}\Big]
\\&
+2{\mathop{\rm 
div}}\Big[\big(\widetilde\mu(a^2)-\widetilde\mu(a^1)\big)
{\mathcal{M}}^1\Big]
+2a^1{\mathop{\rm 
div}}\Big[\big(\widetilde\mu(a^2)-\widetilde\mu(a^1)\big)
{\mathcal{M}}^2\Big]-\frac{1}{2}\delta a\nabla{({B^2})^2}
\\&
-\frac{1}{2}(1+a^1)\nabla((B^2)^2-(B^1)^2)
+(1+a^1)\big(B^2\cdot\nabla\delta B+\delta B\cdot\nabla B^1\big)
+\delta a B^2\cdot\nabla B^2
\end{aligned}
$$
and
$$
\begin{aligned}
G(a^i,u^i,B^i)
&=
B^2\cdot\nabla\delta u+\delta B\cdot\nabla u^1-\delta u\cdot\nabla B^1
+{\mathop{\rm 
div}}\Big\{\big(\widetilde{\sigma}(a^2)-\widetilde{\sigma}(a^1)\big)\nabla B^2\Big\}
\\&
+{\mathop{\rm 
div}}\Big\{\big(\widetilde{\sigma}(a^1)-\sigma^1\big)\nabla\delta B\Big\}.
\end{aligned}
$$
In our discussion we will distinguish between  two cases: the first 
case deals with the situation where 
\hbox{$\frac{1}{p_1}+\frac{1}{p_2}>\frac{2}{N}$}
and the second case concerns
$\frac{1}{p_1}+\frac{1}{p_2}=\frac{2}{N}.$ The distinction between
the two cases appears on the level of the product laws that we use.
\subsection*{\underline{\it The case where $N\geq3,$ $1\leq p_2<2N$ 
and $\frac{1}{p_1}+\frac{1}{p_2}>\frac{2}{N}$}}
We have established the following result.
\begin{prop}\label{CasI}
Let $(a^i,u^i,\nabla\Pi^{i},B^i),$ with $i\in\{ 1,2\},$ be two
solutions of $(\rm\widetilde{MHD})$ system, corresponding to the same initial data 
$a_0\in\dot B^{\frac{N}{p_1}}_{p_1\,\infty}\cap L^\infty(\R^N),$
$u_0,\,B_0\in\dot B^{\frac{N}{p_2}-1}_{p_2\,1}$ with
${\mathop{\rm div}}\, u_0={\mathop{\rm div}}\,B_0=0$ and 
the external forcing term $f$ belonging to
$L^1_{loc}([0,T^*);\,\dot B^{\frac{N}{p_2}-1}_{p_2,1})$ such
that ${\mathcal{Q}}f$ belongs to 
$L^1_{loc}([0,T^{\star});\, \dot B^{\frac{N}{p_2}-2}_{p_2,1}).$ 
Assume that for $i=1,2$ we
have
$$
\begin{aligned}
a^{i}&\in C([0,T^{\star});\,\dot B^{\frac{N}{p_1}}_{p_1,1}(\R^N)),\\
u^{i}&\in C([0,T^{\star});\,\dot B^{\frac{N}{p_2}-1}_{p_2,1})
\cap L^1_{Loc}([0,T^{\star});\,\dot B^{\frac{N}{p}+1}_{p_2,1}),\\
B^{i}&\in C([0,T^{\star});\,\dot B^{\frac{N}{p_2}-1}_{p_2,1})
\cap L^1_{Loc}([0,T^{\star});\,\dot B^{\frac{N}{p}+1}_{p_2,1}),\\
\nabla\Pi^{i}&\in L^1_{Loc}([0,T^{\star});\,\dot 
B^{\frac{N}{p_2}-1}_{p_2,1}).
\end{aligned}
$$
\noindent There exists a positive constant $c$ such that if we have
$$
\Vert a^1\Vert_{L^{\infty}_{T^{\star}}
(\dot B^{\frac{N}{p_1}}_{p_1,\infty}\cap L^\infty)}
\leq
c,
$$
then $(a^2,u^2,\nabla\Pi^2,B^2)=(a^1,u^1,\nabla\Pi^1,B^1).$
\end{prop}
\begin{proof}
The first step of the proof consists in
proving that $(\delta a,\delta u,\nabla\delta\Pi,\delta B) \in
F^p_T,$ where
$$
\begin{aligned}
F^p_T
:&=
C\big([0,T];\,\dot B^{\frac{N}{p_1}-1}_{p_1,1}\big)
\times\big(L^1_T\big(\dot B^{\frac{N}{p_2}}_{p_2,1}\big)
\cap C\big([0,T];\,\dot B^{\frac{N}{p_2}-2}_{p_2,1}\big)\big)
\times\big(L^1_T\big(\dot B^{\frac{N}{p_2}-2}_{p_2,1}\big)\big)
\\&
\times L^1_T\big(\dot B^{\frac{N}{p_2}}_{p_2,1}\big)
\cap C\big([0,T];\,\dot B^{\frac{N}{p_2}-2}_{p_2,1}\big).
\end{aligned}
$$
We define for all $t\leq T$ the quantity
$$
\begin{aligned}
\gamma(t):=\Vert(\delta a,\delta u,\nabla\delta\Pi,&\delta B)
\Vert_{F^p_t}
=
\Vert\delta a\Vert_{L^\infty_t(\dot B^{\frac{N}{p_1}-1}_{p_1,1})}
+\Vert \delta u\Vert_{L^\infty_t(\dot B^{\frac{N}{p_2}-2}_{p_2,1})}
+\mu^1\Vert\delta  u\Vert_{L^1_t(\dot B^{\frac{N}{p_2}}_{p_2,1})}
\\&
+\Vert\nabla\delta \Pi\Vert_{L^1_t(\dot B^{\frac{N}{p_2}-2}_{p_2,1})}
+\Vert \delta B\Vert_{L^\infty_t(\dot B^{\frac{N}{p_2}-2}_{p_2,1})}
+\sigma^1\Vert \delta B\Vert_{L^1_t(\dot B^{\frac{N}{p_2}}_{p_2,1})}.
\end{aligned}
$$
In order to prove that the solution belongs to the space $F_{T}^p,$
it suffices to have 
$(a^{i}-a_0,\overline u^{i},\nabla\overline\Pi^{i},\overline B^{i}) 
\in F^p_T,$ where we
have defined $(\overline u^{i},\nabla\overline\Pi^{i},\overline
B^{i})$ by $ u^{i}=u_L+\overline u^{i},$
$\nabla\Pi^{i}=\nabla\Pi_L^{i} +\nabla\overline\Pi^{i}$ et
$B^{i}=B_L+\overline B^{i}.$ The quantities $u_L,\,\nabla\Pi_L$ and
$B_L$ are defined by the system given bellow:
$$
\left\{\begin{array}{rl}
&\partial_t u_L -\mu^1\Delta u_L +\nabla\Pi_L=f
\medskip\\
&\partial_t B_L -\sigma^1\Delta B_L=0
\medskip \\
&{\mathop{\rm div}}\, u_L ={\mathop{\rm div}}\, B_L=0
\medskip\\
&{(u_L,B_L)}_{|t=0}=(u_0,B_0).
\end{array}
\right.
$$
Indeed, we have by Proposition 2.1 of \cite{CH} that $u_L$ and
$B_L$ have their components in the space $C([0,T];\,\dot
B^{\frac{N}{p_2}-1}_{p_2,1}) \cap L^1(0, T;\,\dot
B^{\frac{N}{p_2}+1}_{p_2,1})$ and 
$\nabla\Pi_L\in L^1(0, T;\,\dot B^{\frac{N}{p_2}-1}_{p_2,1}).$ The quantities $(\overline u^{i},
\nabla\overline\Pi^{i},\overline B^i)$ verify
$$
{\rm(MHD_{mod})}
\left\{
\begin{array}{rl}
&\partial_t\overline u^{i} -\mu^1\Delta \overline u^{i}
+ \nabla\overline\Pi^{i}=K(a^{i},u^{i},\nabla\Pi^{i},B^i)
\medskip\\
&\partial_t\overline B^{i} -\sigma^1\Delta \overline B^{i}=L(u^{i},B^i)
\medskip\\
&{\mathop{\rm div}}\,\overline u^{i}={\mathop{\rm div}}\,\overline
B^{i} =0
\medskip\\
&({\overline u}^{i},{\overline B}^{i})_{|t=0}=(0,0),
\end{array}
\right.
$$
where
$$
\begin{aligned}
K(a^{i},u^{i},\nabla\Pi^{i},B^i)
&=
-u^{i}\cdot\nabla u^{i}+a^{i}\big(\mu^1\Delta u^{i}-\nabla\Pi^{i}\big)
+(1+a^{i}){\mathop{\rm div}}\Big[\big(\widetilde\mu(a^{i})-\mu^1\big)
{\mathcal{M}}^i\Big]
\\&
-\frac{1}{2}(1+a^i)\nabla{{B^i}^2}+(1+a^i)B^i\cdot\nabla B^i
\end{aligned}
$$
and
$$
\begin{aligned}
L(u^i,B^i)=B^i\cdot\nabla u^i-u^i\cdot\nabla B^i
+{\mathop{\rm 
div}}\Big\{\big(\widetilde{\sigma}(a^i)-\sigma^1\big)\nabla B^i\Big\}.
\end{aligned}
$$
We apply the operator $ \mathcal{P}$ to the first equation of the system
$(\rm MHD_{mod})$ and we obtain
\begin{equation}\label{Op_Leray}
\partial_t\overline u^{i} -\mu^1\Delta \overline u^{i}
={\mathcal{P}}\Big(K(a^{i},u^{i},\nabla\Pi^{i},B^i)\Big).
\end{equation}
In the same manner, the divergence operator applied to the same
equation gives
\begin{equation}\label{Ter_Pression}
\begin{aligned}
{\mathop{\rm div}}\big((1+a^i)\nabla\Pi^{i}\big)={\mathop{\rm 
div}}\,{\mathcal{Q}}f
&-{\mathop{\rm div}}\,\Big(u^{i}\cdot\nabla u^{i}+
\frac{1}{2}(1+a^i)\nabla{{B^i}^2}
-(1+a^i)B^i\cdot\nabla B^i\Big)\\
&+
{\mathop{\rm div}}\,\Big(\mu^1(a^{i}\Delta u^{i})+
(1+a^{i}){\mathop{\rm 
div}}\Big[\big(\widetilde\mu(a^{i})-\mu^1\big){\mathcal{M}}^i\Big]\Big),
\end{aligned}
\end{equation}
Combining the inequality (\ref{Interpo}) together with the
hypothesis concerning the solutions stated at the beginning, we find
$u^{i},\,B^i\in L^2_T(\dot B^{\frac{N}{p_2}}_{p_2,1}).$ On the other
hand the inequality (\ref{produit3}) gives that
$u^i\otimes u^i,\;B^i\otimes B^i\;\mbox{and}\; {B^i}^2\in L^2_T(\dot 
B^{\frac{N}{p_2}-1}_{p_2,1})$
for $p_2<2N,$ $N\geq3$ and 
$\frac{1}{p_1}+\frac{1}{p_2}>\frac{2}{N}.$  
Inequality
(\ref{p_1p_2}) then implies that
$$
u^{i}\cdot\nabla u^{i},\;\;a^{i}\Delta u^{i}\in 
L^2_T(\dot B^{\frac{N}{p_2}-2}_{p_2,1}),\;\;
(1+a^i)B^i\cdot\nabla B^i\quad\hbox{and}
\quad (1+a^i)\nabla{B^i}^2
\in
L^2_T(\dot B^{\frac{N}{p_2}-2}_{p_2,1}).
$$
Now the inequality (\ref{p_1p_2}) and Taylor's formula with a remainder 
in the integral form imply
$$
\begin{aligned}
\Big\Vert(1+a^{i}){\mathop{\rm 
div}}\Big[\big(\widetilde\mu(a^{i})-\mu^1\big)
{\mathcal{M}}^i\Big]\Big\Vert_{L^2_T(\dot B^{\frac{N}{p_2}-2}_{p_2,1})}
&\lesssim
\big(1+\Vert a^{i}\Vert_{L^\infty_T(\dot B^{\frac{N}{p_1}}_{p_1,\infty}\cap L^\infty)}\big)
\\&
\times
\Big\Vert\big(\widetilde\mu(a^{i})-\mu^1\big)
{\mathcal{M}}^i\Big\Vert_{L^2_T(\dot B^{\frac{N}{p_2}-1}_{p_2,1})}
\\&
\lesssim
\Vert a^{i}\Vert_{L^\infty_T(\dot B^{\frac{N}{p_1}}_{p_1,\infty}\cap 
L^\infty)}
\Vert u^{i}\Vert_{L^2_T(\dot B^{\frac{N}{p_2}}_{p_2,1})}.
\end{aligned}
$$
We conclude also that the left-hand side term of equality 
(\ref{Ter_Pression})
belongs to $L^2_T(\dot B^{\frac{N}{p_2}-3}_{p_2,1}).$ On the other
hand,  inequality (\ref{p_1p_2}) gives
$$
\|a^{i}\nabla\Pi^{i}\|_{L^2_{T}(\dot B^{\frac{N}{p_2}-2}_{p_2,1})}
\leq
\Vert a^{i}\Vert_{L^\infty_T(\dot B^{\frac{N}{p_1}}_{p_1,\infty}\cap 
L^\infty)}
\|\nabla\Pi^{i}\|_{L^2_{T}(\dot B^{\frac{N}{p_2}-2}_{p_2,1})}.
$$
Consequently, the smallness condition on $a^{i}$ together with 
(\ref{Ter_Pression}) give that
$\nabla\Pi^{i}\in L^2_{T}(\dot B^{\frac{N}{p_2}-2}_{p_2,1}).$ This
allows us to obtain, using the hypothesis concerning $a^{i}$ and the
inequality (\ref{p_1p_2}), that $a^{i}\,\nabla\Pi^{i}\in L^1_T(\dot 
B^{\frac{N}{p_2}-2}_{p_2,1}).$
So we conclude that
$K(a^{i},u^{i},\nabla\Pi^{i},B^i)$ belongs to $L^1_T(\dot 
B^{\frac{N}{p_2}-2}_{p_2,1}).$
In the similar manner we have
$L(u^i,B^i)\in L^1_T(\dot B^{\frac{N}{p_2}-2}_{p_2,1}).$
Since the operator ${\mathcal{P}}$ is continuous on the spaces 
$\dot B^s_{p,r},$
 the terms at the left-hand side of equality (\ref{Op_Leray}) belong to
$L^1_T(\dot B^{\frac{N}{p_2}-2}_{p_2,1}).$ Consequently, applying 
Proposition 2.1 of \cite{CH}, we obtain that 
$\overline u^{i},\overline B^i\in
L^1_{T}(\dot B^{\frac{N}{p_2}}_{p_2,1})\cap C([0,T];\, \dot 
B^{\frac{N}{p_2}-2}_{p_2,1})$
and $\nabla\overline\Pi^{i}\in L^1_{T}(\dot 
B^{\frac{N}{p_2}-2}_{p_2,1}).$
For  $a^i$, we write $\partial_t a^{i}=\,-u^{i}\cdot\nabla a^{i}.$ 
 Since ${1\over p_1}\leq {1\over N}+{1\over p_2}$ 
therefore, the product
laws (\ref{p_1p_2}) allow us to see that $\partial_t a^{i}$ belongs 
to $L^2_{T}(\dot B^{\frac{N}{p_1}-1}_{p_1,1}),$ 
which gives by the inequality of Cauchy-Schwarz  that 
$(a^{i}-a_0)\in C([0,T];\, \dot B^{\frac{N}{p_1}-1}_{p_1,\infty}).$ 
Finally we have $ (\delta a,\delta u,\nabla\delta\Pi,\delta B)\in F^p_T.$
\\


Using these Propositions \ref{eqtransport} and \ref{eq.stokes}) we prove successively that for all $t\leq T$
$$
\Vert\delta a\Vert_{L^\infty_t(\dot B^{\frac{N}{p_1}-1}_{p_1,1})}
\lesssim
e^{C\Vert u^2\Vert_{L^1_t(\dot B^{\frac{N}{p_2}+1}_{p_2,1})}}
\Vert\delta u\Vert_{L^1_t(\dot B^{\frac{N}{p_2}}_{p_2,1})}
\Vert\nabla a^1\Vert_{L^\infty_t(\dot B^{\frac{N}{p_1}-1}_{p_1,1})},
$$
$$
\begin{aligned}
\Vert\delta u\Vert_{L^\infty_t(\dot B^{\frac{N}{p_2}-2}_{p_2,1})}
+\mu^1\Vert\delta u\Vert_{L^1_t(\dot B^{\frac{N}{p_2}}_{p_2,1})}
+\Vert\nabla\delta\Pi&\Vert_{L^1_t(\dot B^{\frac{N}{p_2}-2}_{p_2,1})}
\lesssim
e^{C\Vert u^2\Vert_{L^1_t(\dot B^{\frac{N}{p_2}+1}_{p_2,1})}}
\\&
\times
\Vert H(a^{i},u^{i},\nabla\Pi^{i},B^i)
\Vert_{L^1_t(\dot B^{\frac{N}{p_2}-2}_{p_2,1})}
\end{aligned}
$$
and
$$
\begin{aligned}
\Vert\delta B\Vert_{L^\infty_t(\dot B^{\frac{N}{p_2}-2}_{p_2,1})}
+\sigma^1\Vert\delta B\Vert_{L^1_t(\dot B^{\frac{N}{p_2}}_{p_2,1})}
\lesssim
e^{C\Vert u^2\Vert_{L^1_t(\dot B^{\frac{N}{p_2}+1}_{p_2,1})}}
\Vert G(a^{i},u^{i},B^i)\Vert_{L^1_t(\dot
B^{\frac{N}{p_2}-2}_{p_2,1})}.
\end{aligned}
$$
We will estimate next the term
$H(a^{i},u^{i},\nabla\Pi^{i},B^{i}).$ Inequalities (\ref{produit3})
and (\ref{p_1p_2}) give
$$
\begin{aligned}
\Big\Vert-\delta u&\cdot\nabla u^1+a^1(\mu^1\Delta\delta u-\nabla\delta\Pi)
+\delta a(\mu^1\Delta u^2-\nabla\Pi^2)\Big\Vert_{L^1_T(\dot
B^{\frac{N}{p_2}-2}_{p_2,1})}
\lesssim
\Vert\delta u\Vert_{L^2_T(\dot B^{\frac{N}{p_2}-1}_{p_2,1})}
\\&
\times
\Vert u^1\Vert_{L^2_T(\dot B^{\frac{N}{p_2}}_{p_2,1})}
+\Vert a^1\Vert_{L^\infty_T(\dot B^{\frac{N}{p_1}}_{p_1,\infty}
\cap L^\infty)}
\Big(\Vert\Delta\delta u\Vert_{L^1_T(\dot B^{\frac{N}{p_2}-2}_{p_2,1})}
+\Vert\nabla\delta\Pi\Vert_{L^1_T(\dot B^{\frac{N}{p_2}-2}_{p_2,1})}\Big)
\\&
+\Vert\delta a\Vert_{L^\infty_T(\dot B^{\frac{N}{p_1}-1}_{p_1,\infty})}
\Big(\Vert\Delta u^2\Vert_{L^1_T(\dot B^{\frac{N}{p_2}-1}_{p_2,1})}
+\Vert\nabla\Pi^2\Vert_{L^1_T(\dot B^{\frac{N}{p_2}-1}_{p_2,1})}\Big).
\end{aligned}
$$
Owing to (\ref{p_1p_2}) and Taylor's formula with a remainder in the 
integral form, one finds
\begin{equation}\label{TAYLOR}
\begin{aligned}
\Big\Vert{\mathop{\rm div}}\Big[\big(\widetilde\mu(a^1)-\mu^1\big)
\delta{\mathcal{M}}\Big]+ a^1{\mathop{\rm div}}\Big[\big(\widetilde\mu(a^1)&-\mu^1\big)
\delta{\mathcal{M}}\Big]\Big\Vert_{L^1_T(\dot B^{\frac{N}{p_2}-2}_{p_2,1})}
\\&
\lesssim
\Vert a^1\Vert_{L^\infty_T(\dot B^{\frac{N}{p_1}}_{p_1,\infty}
\cap L^\infty)}
\Vert\delta u\Vert_{L^1_T(\dot B^{\frac{N}{p_2}}_{p_2,1})}
\end{aligned}
\end{equation}
for $p_1\leq p_2.$ Using once more the inequality (\ref{p_1p_2}), 
Taylor's formula, inequality (\ref{produit3}),
and the fact that the space of Besov is stable by the action of a 
$C^\infty$-function
(see for example \cite{M}), one obtains
$$
\begin{aligned}
\Big\Vert{\mathop{\rm div}}\big[\big(\widetilde\mu(a^2)-\widetilde\mu(a^1)\big)
{\mathcal{M}}^2\Big]\Big\Vert_{L^1_{T}
(\dot B^{\frac{N}{p_2}-2}_{p_2,1})}
&\lesssim
\int_0^T\Vert\widetilde\mu(a^2)-\widetilde\mu(a^1)
\Vert_{\dot B^{\frac{N}{p_1}-1}_{p_1,1}}
\Vert\nabla u^2\Vert_{\dot B^{\frac{N}{p_2}}_{p_2,1}}dt
\\&
\lesssim
\int_0^T\Vert\delta a\Vert_{\dot B^{\frac{N}{p_1}-1}_{p_1,1}}
\Vert u^2\Vert_{\dot B^{\frac{N}{p_2}+1}_{p_2,1}}dt.
\end{aligned}
$$
Combining the  inequality (\ref{p_1p_2}) together with an
interpolation result in the temporal variable, we prove that
$$
\begin{aligned}
\Big\Vert\delta a\nabla {(B^2)^2}\Big\Vert_{L^1_T(\dot 
B^{\frac{N}{p_2}-2}_{p_2,1})}
&\lesssim
\Vert\delta a\Vert_{L^{\infty}_T(\dot 
B^{\frac{N}{p_1}-1}_{p_1,\infty})}
\Vert (B^2)^2\Vert_{L^1_T(\dot B^{\frac{N}{p_2}}_{p_2,1})}
\\&
\lesssim
\Vert\delta a\Vert_{L^\infty_T(\dot B^{\frac{N}{p_1}-1}_{p_1,\infty})}
\Vert B^2\Vert_{L^2_T(\dot B^{\frac{N}{p_2}}_{p_2,1})}^2
\\&
\lesssim
\Vert\delta a\Vert_{L^\infty_T(\dot B^{\frac{N}{p_1}-1}_{p_1,\infty})}
\Vert B^2\Vert_{L^\infty_T(\dot B^{\frac{N}{p_2}-1}_{p_2,1})}
\Vert B^2\Vert_{L^1_T(\dot B^{\frac{N}{p_2}+1}_{p_2,1})}.
\end{aligned}
$$
In the same manner we find
$$
\begin{aligned}
\Big\Vert(1+a^1)\nabla((B^2)^2-(B^1)^2)&
\Big\Vert_{L^1_T(\dot B^{\frac{N}{p_2}-2}_{p_2,1})}
\lesssim
\Big(\Vert B^1\Vert_{L^2_T(\dot B^{\frac{N}{p_2}}_{p_2,1})}
+\Vert B^2\Vert_{L^2_T(\dot B^{\frac{N}{p_2}+1}_{p_2,1})}\Big)
\Vert\delta B\Vert_{L^2_T(\dot B^{\frac{N}{p_2}-1}_{p_2,1})}
\\&
\lesssim
\sum_{i=1}^2
\Vert B^i\Vert_{L^2_T(\dot B^{\frac{N}{p_2}}_{p_2,1})}
\Big(\Vert\delta B\Vert_{L^\infty_T(\dot B^{\frac{N}{p_2}-2}_{p_2,1})}
+\Vert\delta B\Vert_{L^1_T(\dot B^{\frac{N}{p_2}}_{p_2,1})}\Big).
\end{aligned}
$$
We have $\frac{1}{p_1}+\frac{1}{p_2}>\frac{2}{N}, $ $p_2<2N$ and
$\vert \frac{1}{p_1}-\frac{1}{p_2}\vert\leq\frac{1}{N},$ so the
inequalities (\ref{p_1p_2}) and (\ref{produit3}) imply
$$
\Vert\delta a\,B^2\cdot\nabla B^2\Vert_{L^1_T
(\dot B^{\frac{N}{p_2}-2}_{p_2,1})}
\lesssim
\Vert\delta a\Vert_{L^\infty_T(\dot B^{\frac{N}{p_1}-1}_{p_1,\infty})}
\Vert B^2\Vert_{L^\infty_T(\dot B^{\frac{N}{p_2}-1}_{p_2,1})}
\Vert B^2\Vert_{L^1_T(\dot B^{\frac{N}{p_2}+1}_{p_2,1})}.
$$
Since one has $\vert\frac{1}{p_1}-\frac{1}{p_2}\vert<\frac{1}{N},$ 
$\frac{1}{p_1}+\frac{1}{p_2}>\frac{2}{N}$
and $p_1\leq p_2,$  Lemmas 3.1 and 3.2 of \cite{HA} remain valid. Thus, combining the preceding inequalities with these Lemmas, we find
$$
\begin{aligned}
\Vert H(a^{i},u^{i},\nabla\Pi^{i},&B^i)\Vert_{L^1_{T}
(\dot B^{\frac{N}{p_2}-2}_{p_2,1})}
\lesssim
 \gamma(t)\Big\{
\Vert (u^1,u^2)\Vert_{L^1_{T}(\dot B^{\frac{N}{p_2}+1}_{p_2,1})}
+\Vert\nabla\Pi^2\Vert_{L^1_{T}(\dot B^{\frac{N}{p_2}-1}_{p_2,1})}
\\&
+
\Vert (B^1,B^2)\Vert_{L^1_T(\dot B^{\frac{N}{p_2}+1}_{p_2,1})}
+\Vert(B^1, B^2)\Vert_{L^2_T(\dot B^{\frac{N}{p_2}+1}_{p_2,1})}
+\Vert a^1\Vert_{L^\infty_{T}(\dot B^{\frac{N}{p_1}}_{p_1,\infty}\cap 
L^\infty)}\Big\}
\\&
+\int_0^{T}\Vert\delta a(t)\Vert_{\dot B^{\frac{N}{p_1}-1}_{p_1,1}}
\Vert u^2\Vert_{\dot B^{\frac{N}{p_2}+1}_{p_2,1}}dt.
\end{aligned}
$$
We need now to estimate $G(a^{i}, u^{i}, B^{i}).$ Since 
${\mathop{\rm div}}\,B^2=0,$ then using the
inequalities of Bernstein and (\ref{produit3}) together with an 
interpolation argument we obtain
$$
\begin{aligned}
\Vert B^2\cdot\nabla\delta u
\Vert_{L^1_T(\dot B^{\frac{N}{p_2}-2}_{p_2,1})}
&\lesssim
\Vert B^2\otimes\delta u\Vert_{L^1_T(\dot B^{\frac{N}{p_2}-1}_{p_2,1})}
\\&
\lesssim
\Vert B^2\Vert_{L^2_T(\dot B^{\frac{N}{p_2}}_{p_2,1})}
\Vert\delta u\Vert_{L^2_T(\dot B^{\frac{N}{p_2}-1}_{p_2,1})}
\\&
\lesssim
\Vert B^2\Vert_{L^2_T(\dot B^{\frac{N}{p_2}}_{p_2,1})}
\Big(\Vert\delta u\Vert_{L^\infty_T(\dot B^{\frac{N}{p_2}-2}_{p_2,1})}
+\Vert\delta u\Vert_{L^1_T(\dot B^{\frac{N}{p_2}}_{p_2,1})}\Big).
\end{aligned}
$$
In the same manner, we have
$$
\begin{aligned}
\Vert\delta B\cdot\nabla u^1-\delta u\cdot\nabla B^1
\Vert_{L^1_T(\dot B^{\frac{N}{p_2}-2}_{p_2,1})}
&\lesssim
\Vert u^1\Vert_{L^2_T(\dot B^{\frac{N}{p_2}}_{p_2,1})}
\Big(\Vert\delta B\Vert_{L^\infty_T(\dot B^{\frac{N}{p_2}-2}_{p_2,1})}
+\Vert\delta B\Vert_{L^1_T(\dot B^{\frac{N}{p_2}}_{p_2,1})}\Big)
\\&
+
\Vert B^1\Vert_{L^2_T(\dot B^{\frac{N}{p_2}}_{p_2,1})}
\Big(\Vert\delta u\Vert_{L^\infty_T(\dot B^{\frac{N}{p_2}-2}_{p_2,1})}
+\Vert\delta u\Vert_{L^1_T(\dot B^{\frac{N}{p_2}}_{p_2,1})}\Big).
\end{aligned}
$$
Arguing similarly  to the case of  inequality (\ref{TAYLOR}), one finds 
that
$$
\Big\Vert{\mathop{\rm div}}
\Big\{\big(\widetilde{\sigma}(a^1)-\sigma^1\big)
\nabla\delta B\Big\}\Big\Vert_{L^1_T(\dot B^{\frac{N}{p_2}-2}_{p_2,1})}
\lesssim
\Vert a^1\Vert_{L^\infty_T(\dot B^{\frac{N}{p_1}}_{p_1,\infty}\cap 
L^\infty)}
\Vert\delta B\Vert_{L^1_T(\dot B^{\frac{N}{p_2}}_{p_2,1})}.
$$
Using the above estimates and Lemmas 3.1, 3.2 of \cite{HA}, and arguing
in the same manner as for the $H$ term, we obtain finally that
$$
\begin{aligned}
\Vert G(a^{i},u^{i},\nabla\Pi^{i},&B^i)\Vert_{L^1_{T}(\dot
B^{\frac{N}{p_2}-2}_{p_2,1})}
\lesssim
\gamma(t)
\bigg(
\Vert (u^1,u^2)\Vert_{L^1_{T}
(\dot B^{\frac{N}{p_2}+1}_{p_2,1})\cap L^2_{T}
(\dot B^{\frac{N}{p_2}}_{p_2,1})}
+
\Vert\nabla\Pi^2\Vert_{L^1_{T}(\dot B^{\frac{N}{p_2}-1}_{p_2,1})}
\\&
+
\Vert (B^1,B^2)\Vert_{L^1_T(\dot B^{\frac{N}{p_2}+1}_{p_2,1})
\cap L^2_{T}(\dot B^{\frac{N}{p_2}}_{p_2,1})}
+
\Vert a^1\Vert_{L^\infty_{T}(\dot B^{\frac{N}{p_1}}_{p_1,\infty}\cap 
L^\infty)}\bigg)
\\&
+\int_0^{T}\|\delta a(t)\|_{\dot B^{{N\over p_1}-1}_{p_1,1}}
\Vert B^2\Vert_{\dot B^{\frac{N}{p_2}+1}_{p_2,1}}dt.
\end{aligned}
$$
Thus, one finds for $t\leq T$ that
$$
\begin{aligned}
\gamma(t)
&\lesssim
\gamma(t)
\bigg(
\Vert (u^1,u^2)\Vert_{L^1_{T}
(\dot B^{\frac{N}{p_2}+1}_{p_2,1})
\cap L^2_{T}(\dot B^{\frac{N}{p_2}}_{p_2,1})}
+\Vert\nabla\Pi^2\Vert_{L^1_{T}(\dot B^{\frac{N}{p_2}-1}_{p_2,1})}
+
\Vert a^1\Vert_{L^\infty_{T}(\dot B^{\frac{N}{p_1}}_{p_1,\infty}\cap 
L^\infty)}
\\&
+
\Vert (B^1,B^2)\Vert_{L^1_T(\dot B^{\frac{N}{p_2}+1}_{p_2,1})
\cap L^2_{T}(\dot B^{\frac{N}{p_2}}_{p_2,1})}
\bigg)
+\int_0^{T} \gamma(t)
\Vert(u^2,B^2)\Vert_{\dot B^{\frac{N}{p_2}+1}_{p_2,1}}dt.
\end{aligned}
$$
We choose a small time $T_{1}\leq T$ such that we have for a constant
$c>0$ small enough the following inequality
$$
\begin{aligned}
\Vert (u^1,u^2)&\Vert_{L^1_{T_{1}}
(\dot B^{\frac{N}{p_2}+1}_{p_2,1})
\cap L^2_{T}(\dot B^{\frac{N}{p_2}}_{p_2,1})}
+\Vert\nabla\Pi^2\Vert_{L^1_{T_{1}}
(\dot B^{\frac{N}{p_2}-1}_{p_2,1})}
\leq
c
\\&
\mbox{and}\quad \Vert (B^1,B^2)\Vert_{L^1_{T_{1}}
(\dot B^{\frac{N}{p_2}+1}_{p_2,1})
\cap L^2_{T}(\dot B^{\frac{N}{p_2}}_{p_2,1})}
\leq
c.
\end{aligned}
$$
Using the assumption that
$\Vert a^1\Vert_{L^\infty_{T_{1}}(\dot 
B^{\frac{N}{p_1}}_{p_1,\infty}\cap L^\infty)}\leq c, $
we have  $\forall t\leq T_{1}$
$$
\gamma(t)\leq
C\int_0^{t} \gamma(t)
\Vert (u^2,B^2)\Vert_{\dot B^{\frac{N}{p_2}+1}_{p_2,1}}dt.
$$
Since the function
 $t\mapsto\Vert u^2\Vert_{\dot B^{\frac{N}{p_2}+1}_{p_2,1}}
 +\Vert B^2\Vert_{\dot B^{\frac{N}{p_2}+1}_{p_2,1}}$
is locally integrable, we deduce by Lemma \ref{osgood}
that $\gamma(t)=0$ for all $t\in[0,T_1]$. It
is easy to see that this property is conserved on the whole time
interval and we obtain finally that $\gamma(t)=0$ for all
$t\in[0,T].$ Thus the proof is complete in the case $1< p_2<2N.$ The
above calculations are  available for $p\neq 1$ (since they are based 
on
Proposition \ref{eq.stokes}). The case $p=1$ is deduced by
injection.
\end{proof}
\subsection*{\underline{\it The case $\frac{1}{p_1}+\frac{1}{p_2}=\frac{2}{N}$ or $N=2$ or $p_2=2N$}}
In this case the condition
$\Vert a^1\Vert_{L^\infty_{T^{\star}}(\dot B^{\frac{N}{p_1}}_{p_1,\infty}\cap L^\infty)}\leq c$
is not  sufficient. To show 
uniqueness, one needs to suppose that
$\Vert a^1\Vert_{L^\infty_{T^{\star}}(\dot B^{\frac{N}{p_1}}_{p_1,1})}
\leq c.$ More precisely, we have the following proposition.
\begin{prop}\label{caslimite}
Let $(a^1,u^1,\nabla\Pi,B^1)$ and $(a^2,u^2,\nabla\Pi^2,B^2)$ be two
solutions of $(\rm\widetilde{MHD})$ corresponding to the initial data
$a_0\in\dot B^{\frac{N}{p_1}}_{p_1,1},$
$u_0,\,B_0\in\dot B^{\frac{N}{p_2}-1}_{p_2,1}$ where 
${\mathop{\rm div}} \,u_0={\mathop{\rm div}}\,B_0=0$ 
and $f$ is such that  its components are in 
$L^1_{loc}([0,T^*);\,\dot B^{\frac{N}{p_2}-1}_{p_2,1})$
and ${\mathcal{Q}}f$ belongs to $L^1_{loc}([0,T^{\star});\,\dot
B^{\frac{N}{p_2}-2}_{p_2,1}).$ We assume that for $i=1,2$ we
have
$$
\begin{aligned}
a^{i}&\in C([0,T^{\star});\,{\mathcal S}')
\cap L^{\infty}_{loc}([0,T^{\star});\,\dot
B^{\frac{N}{p_1}}_{p_1,1}) ,\\
u^{i}&\in C([0,T^{\star});\,\dot B^{\frac{N}{p_2}-1}_{p_2,1})
\cap L^1_{loc}([0,T^{\star});\,\dot
B^{\frac{N}{p_2}+1}_{p_2,1}),\\
B^{i}&\in C([0,T^{\star});\,\dot B^{\frac{N}{p_2}-1}_{p_2,1})
\cap L^1_{loc}([0,T^{\star});\,\dot
B^{\frac{N}{p_2}+1}_{p_2,1}),\\
\nabla\Pi^{i}&\in L^1_{loc}([0,T^{\star});\,\dot
B^{\frac{N}{p_2}}_{p_2,1}).
\end{aligned}
$$
Then there exists a positive constant $c$ which does not depend on 
these
solutions such that the inequality
$$
\Vert a^1\Vert_{\widetilde L^{\infty}_{T^{\star}}
(\dot B^{\frac{N}{p_1}}_{p_1,1})}
\leq
c
$$
implies $(a^2,u^2,\nabla\Pi^2,B^2)=(a^1,u^1,\nabla\Pi^1,B^1).$
\end{prop}
\begin{proof}
We need to prove first that 
$(\delta a,\delta u,\nabla\delta\Pi,\delta B)\in G_T,$ where
$$
G_T\hspace{-,15cm}:=\hspace{-0,15cm}L^\infty_T(\dot B^{\frac{N}{p_1}-1}_{p_1,\infty})
\times
\widetilde{L}^1_T(\dot B^{\frac{N}{p_2}}_{p_2,\infty})
\cap L^\infty_T(\dot B^{-2+\frac{N}{p_2}}_{p_2,\infty})
\times
\widetilde{L}^1_T(\dot B^{-2+\frac{N}{p_2}}_{p_2,\infty})
\times
\widetilde{L}^1_T(\dot B^{\frac{N}{p_2}}_{p_2,\infty})
\cap L^\infty_T(\dot B^{-2+\frac{N}{p_2}}_{p_2,\infty}).
$$
The estimates in this space would allow us to obtain
uniqueness of the solution by the Osgood Lemma. We define
$$
\begin{aligned}
\gamma(t):=
\Vert\delta u\Vert_{L^\infty_T(\dot B^{-2+\frac{N}{p_2}}_{p_2,\infty})}
&+\Vert\delta u\Vert_{\widetilde{L}^1_T(\dot
B^{\frac{N}{p_2}}_{p_2,\infty})}
+\Vert\nabla\delta\Pi\Vert_{\widetilde{L}^1_T(\dot
B^{-2+\frac{N}{p_2}}_{p_2,\infty})}
+\Vert\delta B\Vert_{L^\infty_T(\dot
B^{-2+\frac{N}{p_2}}_{p_2,\infty})}
\\&
+\Vert\delta B\Vert_{\widetilde{L}^1_T
(\dot B^{\frac{N}{p_2}}_{p_2,\infty})}.
\end{aligned}
$$
The term $G_ {T} $ is dealt with in the same way as in the first case. 
The only difference to be noted appears in the treatment of the 
products of the type $a^ {i} \nabla\Pi^{i}. $ Here inequality 
(\ref{sommenulle})
should be used to ensure that the left-hand side term of  equality 
(\ref{Op_Leray}) belongs to $L^2_T(\dot B^ {- 1} _{3,\infty}).$ Thus
Proposition 2.1 of \cite{CH} implies that  $(\delta a, \delta u, 
\nabla\delta\Pi, \delta B)\in G_T.$

In this case it is enough to study the case $\frac{2}{N} =\frac{1}{p_1} 
+\frac{1}{p_2},$ since one can deduce the other cases from this one. 
Indeed,  if $p_2=2N,$ then $p_1= \frac{2N}{3},$ since
$\frac{1}{p_1}\leq\frac{1}{N}+\frac{1}{p_2}$ and 
$\frac{2}{N}\leq\frac{1}{p_1}+\frac{1}{p_2}.$ Therefore it is a particular case of 
$\frac{2}{N}=\frac{1}{p_1}+\frac{1}{p_2}.$ For $N=2,$ one starts with $p_2=4$ and 
$p_1=\frac{4}{3}$. Afterwards by injection, one will have uniqueness 
for $1\leq p_1\leq\frac{4}{3}$ and
$1\leq p_2\leq 4,$ the same for $1\leq p_1\leq 4$ and $1\leq 
p_2\leq\frac{4}{3}. $ Hence one can suppose that $\frac{2}{N}=\frac{1}{p_1} + 
\frac{1}{p_2}.$ Moreover, one can suppose that $p_2\geq2$ since
inequality (\ref{es.st.li}) is valid for $p\geq2$. The case $p_2\leq2$ 
follows by injection.\\
Using Propositions \ref{eqtransport} and \ref{eq.stokes}, we have
\begin{equation}\label{delta_a}
\Vert\delta a\Vert_{L^\infty_t(\dot B^{\frac{N}{p_1}-1}_{p_1,\infty})}
\leq
e^{C\Vert\nabla u^2\Vert_{L^1_t(\dot B^{\frac{N}{p_2}}_{p_2,1})}}
\Vert\delta u\cdot\nabla a^1\Vert_{L^1_t
(\dot B^{\frac{N}{p_1}-1}_{p_1,\infty})},
\end{equation}
\begin{equation}\label{delta_u}
\begin{aligned}
\Vert\delta u\Vert_{L^\infty_t(\dot B^{-2+\frac{N}{p_2}}_{p_2,\infty})}
+\mu^1\Vert\delta u\Vert_{\widetilde{L}^1_t(\dot
B^{\frac{N}{p_2}}_{p_2,\infty})}
&+\Vert\nabla\delta\Pi\Vert_{\widetilde{L}^1_t(\dot
B^{-2+\frac{N}{p_2}}_{p_2,\infty})}
\leq
Ce^{C\Vert\nabla u^2\Vert_{L^1_t(\dot B^{\frac{N}{p_2}}_{p_2,1})}}
\\&
\times
\Vert H(a^i,u^i,\nabla\Pi^i,B^i)\Vert_{\widetilde {L}^1_T
(\dot B^{-2+\frac{N}{p_2}}_{p_2,\infty})}
\end{aligned}
\end{equation}
and
\begin{equation}\label{dell}
\begin{aligned}
\Vert\delta B\Vert_{L^\infty_t(\dot B^{-2+\frac{N}{p_2}}_{p_2,\infty})}
+\sigma^1\Vert\delta u\Vert_{\widetilde{L}^1_t
(\dot B^{\frac{N}{p_2}}_{p_2,\infty})}
&\leq
Ce^{C\Vert\nabla u^2\Vert_{L^1_t(\dot B^{\frac{N}{p_2}}_{p_2,1})}}
\Vert G(a^i,u^i,B^i)\Vert_{\widetilde {L}^1_T
(\dot B^{-2+\frac{N}{p_2}}_{p_2,\infty})}.
\end{aligned}
\end{equation}
Combining the estimates of  $\delta a,$ inequality
(\ref{p_1p_2}), the Bernstein and  Minkowski inequalities, we obtain
\begin{equation}\label{produit_1}
\Vert\delta u\cdot\nabla a^1\Vert_{L^1_t
(\dot B^{\frac{N}{p_1}-1}_{p_1,\infty})}
\lesssim
\Vert\delta u\Vert_{L^1_t(\dot B^{\frac{N}{p_2}}_{p_2,1})}
\Vert a^1\Vert_{\widetilde L^\infty_t(\dot B^{\frac{N}{p_1}}_{p_1,1})}.
\end{equation}
By Lemma \ref{DANCHIN} one has
\begin{equation}\label{produit_2}
\begin{aligned}
\Vert\delta u&\Vert_{L^1_t(\dot B^{\frac{N}{p_2}}_{p_2,1})}
\lesssim
\Vert\delta u\Vert_{\widetilde{L}^1_t
(\dot B^{\frac{N}{p_2}}_{p_2,\infty})}
\log\Big(e+\frac{\Vert\delta u\Vert_{\widetilde{L}^1_t
(\dot B^{\frac{N}{p_2}-1}_{p_2,\infty})}
+
\Vert\delta u\Vert_{\widetilde{L}^1_t
(\dot B^{1+\frac{N}{p_2}}_{p_2,\infty})}}
{\Vert\delta u\Vert_{\widetilde{L}^1_t
(\dot B^{\frac{N}{p_2}}_{p_2,\infty})}}\Big)
\\&
\lesssim
\Vert\delta u\Vert_{\widetilde{L}^1_t
(\dot B^{\frac{N}{p_2}}_{p_2,\infty})}
\log\Big(e+\frac{t\sum_{i=1}^2\Vert u^i\Vert_{L^\infty_t
(\dot B^{\frac{N}{p_2}-1}_{p_2,\infty})}
+
\sum_{i=1}^2\Vert u^i\Vert_{L^1_t(\dot B^{1+\frac{N}{p_2}}_{p_2,1})}}
{\Vert\delta u\Vert_{\widetilde{L}^1_t
(\dot B^{\frac{N}{p_2}}_{p_2,\infty})}}\Big).
\end{aligned}
\end{equation}
We will now  estimate the term $H(a^{i},u^{i},\nabla\Pi^{i},B^{i}).$
Since ${\mathop{\rm div}}\,\delta u=0$ the inequalities of Bernstein, 
(\ref{produit3}) (for $p_2<2N$) and (\ref{sommenulle}) for
($p_2=2N$) imply
$$
\begin{aligned}
\Vert\delta u\cdot\nabla u^1\Vert_{\widetilde L^1_t
(\dot B^{-2+\frac{N}{p_2}}_{p_2,\infty})}
&\lesssim
\Vert\delta u\otimes u^1\Vert_{\widetilde L^1_t
(\dot B^{-1+\frac{N}{p_2}}_{p_2,\infty})}
\\&
\lesssim
\Vert u^1\Vert_{\widetilde L^2_t(\dot B^{\frac{N}{p_2}}_{p_2,1})}
\Vert\delta u\Vert_{\widetilde L^2_t
(\dot B^{-1+\frac{N}{p_2}}_{p_2,\infty})}
\\&
\lesssim
\Vert u^1\Vert_{\widetilde L^2_t(\dot B^{\frac{N}{p_2}}_{p_2,1})}
\big(\Vert\delta u\Vert_{\widetilde L^1_t
(\dot B^{\frac{N}{p_2}}_{p_2,\infty})}
+\Vert\delta u\Vert_{L^\infty_t
(\dot B^{-2+\frac{N}{p_2}}_{p_2,\infty})}\big).
\end{aligned}
$$
Since $p_1\leq p_2,$ $\frac{2}{N}\leq\frac{1}{p_1}+\frac{1}{p_2}\leq1,$  
thanks to inequalities
(\ref{sommenulle}) and using the Bernstein inequality, we have
$$
\begin{aligned}
\Big\Vert
a^1\big(\mu^1\Delta\delta u-\nabla\delta\Pi\big)
&+\delta a\big(\mu^1\Delta u^2-\nabla\Pi^2\big)
\Big\Vert_{\widetilde L^1_t(\dot B^{-2+\frac{N}{p_2}}_{p_2,\infty})}
\lesssim
\Vert a^1\Vert_{\widetilde L^\infty_t(\dot B^{\frac{N}{p_1}}_{p_1,1})}
\\&
\times
\big(\Vert\delta u\Vert_{\widetilde L^1_t
(\dot B^{\frac{N}{p_2}}_{p_2,\infty})}
+\Vert\nabla\delta\Pi\Vert_{\widetilde L^1_t
(\dot B^{-2+\frac{N}{p_2}}_{p_2,\infty})}\big)
\\&
+\int_0^t\Vert\delta a\Vert_{\dot B^{-1+\frac{N}{p_1}}_{p_1,\infty}}
\Big(\Vert u^2\Vert_{\dot B^{1+\frac{N}{p_2}}_{p_2,1}}
+\Vert\nabla\Pi^2\Vert_{\dot B^{-1+\frac{N}{p_2}}_{p_2,1}}\Big)d\tau.
\end{aligned}
$$
The  inequality of Minkowski, (\ref{sommenulle}), (\ref{p_1p_2}), and 
the Taylor formula imply
$$
\begin{aligned}
\Big\Vert\delta a\,{\mathop{\rm
div}}\Big[\big(\widetilde\mu(a^2)-\mu^1\big)
&{\mathcal{M}}^2\Big]\Big\Vert_{\widetilde L^1_t
(\dot B^{-2+\frac{N}{p_2}}_{p_2,\infty})}
\lesssim
\int_0^t\Vert\delta a\Vert_{\dot B^{-1+\frac{N}{p_1}}_{p_1,\infty}}
\Big\Vert\big(\widetilde\mu(a^2)-\mu^1\big){\mathcal{M}}^2
\Big\Vert_{\dot B^{\frac{N}{p_2}}_{p_2,1}}d\tau
\\&
\lesssim
\int_0^t\Vert\delta a\Vert_{\dot B^{-1+\frac{N}{p_1}}_{p_1,\infty}}
\Vert \widetilde\mu(a^2)-\mu^1\Vert_{\dot B^{\frac{N}{p_1}}_{p_1,1}}
\Vert\nabla u^2\Vert_{\dot B^{\frac{N}{p_2}}_{p_2,1}}d\tau
\\&
\lesssim
\Vert a^2\Vert_{\widetilde L^\infty_t(\dot B^{\frac{N}{p_1}}_{p_1,1})}
\int_0^t\Vert\delta a\Vert_{\dot
B^{-1+\frac{N}{p_1}}_{p_1,\infty}}
\Vert u^2\Vert_{\dot B^{1+\frac{N}{p_2}}_{p_2,1}}d\tau.
\end{aligned}
$$
Using the Bernstein inequality and (\ref{p_1p_2}), we find
\begin{equation}\label{deltaB}
\begin{aligned}
\Big\Vert{\mathop{\rm div}}\Big[\big(\widetilde\mu(a^1)-\mu^1\big)
\delta{\mathcal{M}}\Big]\Big\Vert_{\widetilde L^1_t
(\dot B^{-2+\frac{N}{p_2}}_{p_2,\infty})}
&\lesssim
\Vert\widetilde\mu(a^1)-\mu^1\Vert_{\widetilde L^\infty_t
(\dot B^{\frac{N}{p_1}}_{p_1,1})}
\Vert\nabla\delta u\Vert_{\widetilde L^1_t
(\dot B^{-1+\frac{N}{p_2}}_{p_2,\infty})}
\\&
\lesssim
\Vert a^1\Vert_{\widetilde L^\infty_t(\dot B^{\frac{N}{p_1}}_{p_1,1})}
\Vert\delta u\Vert_{\widetilde L^1_t
(\dot B^{\frac{N}{p_2}}_{p_2,\infty})}.
\end{aligned}
\end{equation}
This and the inequality of Minkowski, (\ref{sommenulle}), the Bernstein 
inequality, (\ref{p_1p_2}), Taylor's formula and (\ref{produit3}) give
\begin{equation}\label{a1_div}
\begin{aligned}
\Big\Vert a^1{\mathop{\rm div}}\big[\big(\widetilde\mu(a^2)-\widetilde\mu(a^1)\big)
{\mathcal{M}}^2\big]&\Big\Vert_{\widetilde L^1_t(\dot B^{-2+\frac{N}{p_2}}_{p_2,\infty})}
\\&
\lesssim
\int_0^t\Vert a^1\Vert_{\dot B^{\frac{N}{p_1}}_{p_1,1}}
\Big\Vert\big(\widetilde\mu(a^2)-\widetilde\mu(a^1)\big){\mathcal{M}}^2
\Big\Vert_{\dot B^{-1+\frac{N}{p_2}}_{p_2,\infty}}d\tau
\\&
\lesssim
\int_0^t
\big\Vert\widetilde\mu(a^2)-\widetilde\mu(a^1)
\big\Vert_{\dot B^{-1+\frac{N}{p_1}}_{p_1,\infty}}\Vert\nabla
u^2\Vert_{\dot B^{\frac{N}{p_2}}_{p_2,\infty}
\cap L^\infty}d\tau
\\&
\lesssim
\int_0^t
\Vert\delta a\Vert_{\dot B^{-1+\frac{N}{p_1}}_{p_1,\infty}}
\sum_{i=1}^2\Vert a^i\Vert_{\dot B^{\frac{N}{p_1}}_{p_1,1}}
\Vert u^2\Vert_{ \dot B^{1+\frac{N}{p_2}}_{p_2,1}}d\tau
\\&
\lesssim
\int_0^t\Vert\delta a\Vert_{\dot B^{-1+\frac{N}{p_1}}_{p_1,\infty}}
\Vert u^2\Vert_{ \dot B^{1+\frac{N}{p_2}}_{p_2,1}}d\tau.
\end{aligned}
\end{equation}
In the same manner we obtain the following estimates
$$
\begin{aligned}
\Big\Vert a^1{\mathop{\rm div}}\Big[\big(\widetilde\mu(a^1)-\mu^1\big)
\delta{\mathcal{M}}\Big]\Big\Vert_{\widetilde L^1_t(\dot
B^{-2+\frac{N}{p_2}}_{p_2,\infty})}
&\lesssim
\Vert a^1\Vert_{\widetilde L^\infty_t(\dot B^{\frac{N}{p_1}}_{p_1,1})}
\Big\Vert\big(\widetilde\mu(a^1)-\mu^1\big)
\delta{\mathcal{M}}\Big\Vert_{\widetilde L^1_t(\dot
B^{-1+\frac{N}{p_2}}_{p_2,\infty})}
\\&
\lesssim
\Vert a^1\Vert_{\widetilde L^\infty_t(\dot
B^{\frac{N}{p_1}}_{p_1,1})}^2
\Vert\delta u\Vert_{\widetilde L^1_t(\dot
B^{\frac{N}{p_2}}_{p_2,\infty})}
\end{aligned}
$$
and
\begin{equation}\label{magnetique}
\Big\Vert{\mathop{\rm
div}}\Big[\big(\widetilde\mu(a^2)-\widetilde\mu(a^1)\big)
{\mathcal{M}}^2\Big]\Big\Vert_{\widetilde L^1_t
(\dot B^{-2+\frac{N}{p_2}}_{p_2,\infty})}
\lesssim
\int_0^t\Vert\delta a\Vert_{\dot B^{-1+\frac{N}{p_1}}_{p_1,\infty}}
\Vert u^2\Vert_{\dot B^{1+\frac{N}{p_2}}_{p_2,1}}d\tau.
\end{equation}
Using the Minkowski  inequality,  (\ref{p_1p_2}), the fact that
 $\dot{B}_{p_2,1}^{\frac{N}{p_2}}$ is an algebra, and interpolation, we 
obtain
$$
\begin{aligned}
\Vert\delta a\nabla{(B^2)^2}\Vert_{\widetilde L^1_t
(\dot B^{-2+\frac{N}{p_2}}_{p_2,\infty})}
&\lesssim
\int_{0}^t\Vert\delta a\Vert_{\dot B^{-1+\frac{N}{p_1}}_{p_1,\infty}}
\Vert {B^2}\Vert_{\dot B^{\frac{N}{p_2}}_{p_2,1}}^2d\tau\\
&\lesssim
\int_{0}^t\Vert\delta a\Vert_{\dot B^{-1+\frac{N}{p_1}}_{p_1,\infty}}
\Vert B^2\Vert_{\dot B^{-1+\frac{N}{p_2}}_{p_2,1}}
\Vert B^2\Vert_{\dot B^{1+\frac{N}{p_2}}_{p_2,1}}d\tau.
\end{aligned}
$$
Since ${\mathop{\rm div}}\,B^2=0$, in an analogous manner, we obtain
$$
\begin{aligned}
\Vert\delta a B^2\cdot\nabla B^2\Vert_{\widetilde L^1_t
(\dot B^{-2+\frac{N}{p_2}}_{p_2,\infty})}
&\lesssim\int_{0}^t\Vert\delta a
\Vert_{\dot B^{-1+\frac{N}{p_1}}_{p_1,\infty}}
\Vert B^2\Vert_{\dot B^{-1+\frac{N}{p_2}}_{p_2,1}}
\Vert B^2\Vert_{\dot B^{1+\frac{N}{p_2}}_{p_2,1}}d\tau.
\end{aligned}
$$
Thanks to inequalities (\ref{sommenulle}), (\ref{produit3})
and a classical interpolation argument, we
can write
$$
\begin{aligned}
\Big\Vert\big(1+a^1\big)\nabla\big((B^2)^2-(B^1)^2\big)
&\Big\Vert_{\widetilde L^1_t(\dot B^{-2+\frac{N}{p_2}}_{p_2,\infty})}
\lesssim
\Big(1+\Vert a^1\Vert_{\widetilde L^\infty_t
(\dot B^{\frac{N}{p_1}}_{p_1,1})}\Big)
\Big\Vert(B^2)^2-(B^1)^2\Big\Vert_{\widetilde L^1_t
(\dot B^{-1+\frac{N}{p_2}}_{p_2,\infty})}
\\&
\lesssim
\sum_{i=1}^2\Vert B^i\Vert_{\widetilde L^2_t
(\dot B^{\frac{N}{p_2}}_{p_2,1})}
\Vert\delta B\Vert_{\widetilde L^2_t
(\dot B^{-1+\frac{N}{p_2}}_{p_2,\infty})}
\\&
\lesssim
\sum_{i=1}^2\Vert B^i\Vert_{\widetilde L^2_t
(\dot B^{\frac{N}{p_2}}_{p_2,1})}
\Big(\Vert\delta B\Vert_{L^\infty_t
(\dot B^{-2+\frac{N}{p_2}}_{p_2,\infty})}
+\Vert\delta B\Vert_{\widetilde L^1_t
(\dot B^{\frac{N}{p_2}}_{p_2,\infty})}\Big).
\end{aligned}
$$
Since ${\mathop{\rm div}}\,\delta B={\mathop{\rm div}}\,B^2=0,$ one 
will have in the same way
$$
\begin{aligned}
\Big\Vert\big(1+a^1\big)\big(\delta B\cdot\nabla 
B^1+B^2\cdot\nabla\delta B\big)
\Big\Vert_{\widetilde L^1_t(\dot B^{-2+\frac{N}{p_2}}_{p_2,\infty})}
&\lesssim
\sum_{i=1}^2\Vert B^i\Vert_{\widetilde L^2_t
(\dot B^{\frac{N}{p_2}}_{p_2,1})}
\\&
\times
\Big(\Vert\delta B\Vert_{L^\infty_t
(\dot B^{-2+\frac{N}{p_2}}_{p_2,\infty})}
+\Vert\delta B\Vert_{\widetilde L^1_t
(\dot B^{\frac{N}{p_2}}_{p_2,\infty})}\Big).
\end{aligned}
$$
Combining all these estimates, we are able to establish
\begin{equation}\label{Lemma-1}
\begin{aligned}
&\Vert H(a^{i},u^{i},\nabla\Pi^{i},B^i)\Vert_{\widetilde L^1_t
(\dot B^{-2+\frac{N}{p_2}}_{p_2,\infty})}
\lesssim
\gamma(t)
\bigg(
\Vert(u^1,u^2)\Vert_{L^1_t(\dot B^{1+\frac{N}{p_2}}_{p_2,1})
\cap \widetilde L^2_t(\dot B^{\frac{N}{p_2}}_{p_2,1})}
\\&
+\Vert a^1\Vert_{\widetilde L^\infty_t(\dot B^{\frac{N}{p_1}}_{p_1,1})}
+\Vert (B^1,B^2)\Vert_{L^1_t(\dot B^{1+\frac{N}{p_2}}_{p_2,1})
\cap\widetilde L^2_t(\dot B^{\frac{N}{p_2}}_{p_2,1})}\bigg)
\\&
+\int_{0}^t\|\delta 
a\|_{\dot{B}_{p_1,\infty}^{-1+\frac{N}{p_1}}}g(\tau)d\tau,
\end{aligned}
\end{equation}
where $g$ is a locally integrable function.\\
We give now the estimates for $G.$ Using the Bernstein inequality and
(\ref{produit3}) for $p_2<2N,$ (\ref{sommenulle}) for $p_2=2N$, we 
obtain by interpolation
$$
\begin{aligned}
\Big\Vert B^2\cdot\nabla\delta u&+\delta B\cdot\nabla u^1-\delta
u\cdot\nabla B^1
\Big\Vert_{\widetilde L^1_t(\dot B^{-2+\frac{N}{p_2}}_{p_2,\infty})}
\\&
\lesssim
\Vert B^2\otimes\delta u+\delta B\otimes u^1-\delta u\otimes B^1
\Vert_{\widetilde L^1_t(\dot B^{-1+\frac{N}{p_2}}_{p_2,\infty})}
\\&
\lesssim
\Vert (B^1,B^2)\Vert_{\widetilde L^2_t(\dot B^{\frac{N}{p_2}}_{p_2,1})}
\Vert\delta u\Vert_{\widetilde L^2_t(\dot 
B^{-1+\frac{N}{p_2}}_{p_2,\infty})}
+\Vert\delta B\Vert_{\widetilde L^2_t(\dot 
B^{-1+\frac{N}{p_2}}_{p_2,\infty})}
\Vert u^1\Vert_{\widetilde L^2_t(\dot B^{\frac{N}{p_1}}_{p_2,1})}
\\&
\lesssim
\Vert(B^1,B^2,u^1)\Vert_{\widetilde L^2_t
(\dot B^{\frac{N}{p_2}}_{p_2,1})}
\\&
\times
\Big(\Vert\delta B\Vert_{L^\infty_t(\dot 
B^{-2+\frac{N}{p_2}}_{p_2,\infty})}
+\Vert\delta B\Vert_{\widetilde L^1_t(\dot 
B^{\frac{N}{p_2}}_{p_2,\infty})}
+\Vert\delta u\Vert_{L^\infty_t(\dot 
B^{-2+\frac{N}{p_2}}_{p_2,\infty})}
+\Vert\delta u\Vert_{\widetilde L^1_t(\dot 
B^{\frac{N}{p_2}}_{p_2,\infty})}\Big).
\end{aligned}
$$
We obtain identically to (\ref{magnetique}) and (\ref{deltaB}) that
\begin{eqnarray*}
\Big\Vert{\mathop{\rm 
div}}\Big[\big(\widetilde\sigma(a^2)-\widetilde\sigma(a^1)\big)
\nabla B^2\Big]\Big\Vert_{\widetilde L^1_t
(\dot B^{-2+\frac{N}{p_2}}_{p_2,\infty})}
&\lesssim&
\int_{0}^t\Vert\delta a\Vert_{\dot 
B^{-1+\frac{N}{p_1}}_{p_1,\infty}}\Vert
B^2\Vert_{ \dot B^{1+\frac{N}{p_2}}_{p_2,1}}d\tau
\end{eqnarray*}
and
$$
\Big\Vert{\mathop{\rm 
div}}\Big[\big(\widetilde\sigma(a^1)-\sigma^1\big)
\nabla\delta B\Big]\Big\Vert_{\widetilde L^1_t
(\dot B^{-2+\frac{N}{p_2}}_{p_2,\infty})}
\lesssim
\Vert a^1\Vert_{\widetilde L^\infty_t(\dot B^{\frac{N}{p_1}}_{p_1,1})}
\Vert\delta B\Vert_{\widetilde L^1_t
(\dot B^{\frac{N}{p_2}}_{p_2,\infty})}.
$$
We deduce from these estimates that
$$
\begin{aligned}
\Vert G(a^i,u^i,B^i)\Vert_{\widetilde L^1_t
(\dot B^{-2+\frac{N}{p_2}}_{p_2,\infty})}
&\lesssim
\gamma(t)
\Big(\Vert (u^1,B^1,B^2)\Vert_{\widetilde L^2_t
(\dot B^{\frac{N}{p_2}}_{p_2,1})}
+\Vert a^1\Vert_{\widetilde L^\infty_t(\dot 
B^{\frac{N}{p_1}}_{p_1,1})}\Big)
\\&
+
\int_{0}^t\|\delta a\|_{\dot{B}_{p_1,\infty}^{-1+\frac{N}{p_1}}}
\Vert B^2\Vert_{\dot B^{1+\frac{N}{p_2}}_{p_2,1}}d\tau.
\end{aligned}
$$
Using the above estimate together with those given by 
(\ref{Lemma-1}), 
we have
$$
\begin{aligned}
\gamma(t)
\lesssim
\gamma(t)
\bigg(
\Vert(u^1,u^2,B^1,B^2)\Vert_{L^1_t(\dot B^{1+\frac{N}{p_2}}_{p_2,1})
\cap \widetilde L^2_t(\dot B^{\frac{N}{p_2}}_{p_2,1})}
&+\Vert a^1\Vert_{\widetilde L^\infty_t(\dot B^{\frac{N}{p_1}}_{p_1,1})}\bigg)
\\&
+\int_{0}^t\|\delta a\|_{\dot{B}_{p_1,1}^{-1+\frac{N}{p_1}}}g(\tau)d\tau.
\end{aligned}
$$
Using the above estimate, we may
choose a sufficiently small  time $T_{1}$ so that using inequalities
(\ref{delta_a}), (\ref{produit_1}), (\ref{produit_2}) and the smallness 
of $a^1$, we obtain for all
$t\in[0,T_{1}]$
$$
\gamma(t)\lesssim \int_{0}^t\log\Big(e+\frac{\alpha(T)} {\Vert\delta
u\Vert_{\widetilde{L}^1_{\tau}(\dot
B^{\frac{N}{p_2}}_{p_2,\infty})}}\Big)\Vert\delta u
\Vert_{\widetilde{L}^1_{\tau}(\dot
B^{\frac{N}{p_2}}_{p_2,\infty})}g(\tau)d\tau,$$ with
$\alpha(T)=\sum_{i=1}^2T\|u^i\|_{L^\infty_{T}(\dot{B}_{p_2,1}
^{-1+\frac{N}{p_2}})}+
\|u^i\|_{{L}^1_{T}(\dot{B}_{p_2,1}^{1+\frac{N}{p_2}})}.$ Owing to the 
fact
that \hbox{$x\longmapsto x\log(e+\frac{\alpha(T)}{x})$} is an increasing 
function on
$\R_{+},$ we have for all $t\in[0,T_{1}]$
 $$
\gamma(t)\lesssim \int_{0}^t\gamma(\tau)\log\Big(e+\frac{\alpha(T)}
{\gamma(\tau)}\Big)g(\tau)d\tau.
$$ 
So by Lemma \ref{osgood}, we deduce that
that $\gamma(t)=0,$ for all $t\in[0,T_{1}].$ This gives by
 inequality (\ref{delta_a}), that $\delta a=0$.
Standard arguments now yield the required conclusion. 
We note that the method used in this section (the logarithmic interpolation argument and the application of the Osgood lemma) is inspired by the proofs of the uniqueness given by Danchin \cite{Danchin} and was used by the authors in \cite{Abidi-Paicu}.
\end{proof}
\subsection{Existence.}

\smallskip

\noindent Throughout this section we assume that
$p_1\leq p_2,$ $\frac{1}{p_1}+\frac{1}{p_2}>\frac{1}{N}$ and
$\frac{1}{p_1}\leq \frac{1}{N}+{1\over p_2}.$

\smallskip

\noindent The proof of existence of a solution is performed in a
standard manner. We begin by solving an approximate problem  
and we
prove that the solutions are uniformly bounded. The last step consists
in studying the convergence to a solution of the initial equation.
\subsubsection*{Construction of a regular approximate solution.}
Let us recall first the following result (see [\cite{HA}, Lemma 4.2).
\begin{lem}\label{bravo}{\it
Assume that  $s_{i}\in\R$ and $(p_{i},r_{i})\in[1,\infty[^2$ for
${i=1,2}.$ Let $G\in\dot B^{s_{1}}_{p_{1}\,r_{1}}(\R^N).$ Then there
exists $G^n\in H^\infty(\R^N),$ such that for all $\varepsilon>0$ there
is $n_0$ such that
$$ \Vert G^n- G\Vert_{\dot B^{s_1}_{p_1\,r_1}}
\leq \varepsilon\quad\forall\;n\geq n_0.
$$
If we have ${\mathop{\rm div}}\, G=0$ and 
${\mathcal{Q}}G\in\dot B^{s_2}_{p_2\,r_2},$ then we can choose 
$G^n$ such that ${\mathop{\rm div}}\, G^n=0$ and 
${\mathcal{Q}}G^n$ is uniformly bounded with respect to 
$n$ in the
space $\dot B^{s_2}_{p_2\,r_2}.$}
\end{lem}
Owing to the above Lemma there exist 
$a_0^n ,u_{0}^n,\,B_0^n\in H^\infty(\R^N)$
and $f^n\in L^1_{loc}(\R_+;\,H^\infty(\R^N))$ such that we have
$$
\begin{aligned}
\Vert a^n_0\Vert_{L^\infty}
\lesssim
\Vert a_0\Vert_{L^\infty},
&\quad
{\mathop{\rm div}} \,u^n_0={\mathop{\rm div}}\,B^n_0=0
\\&
\quad\hbox{and}\quad
\Vert{\mathcal{Q}}f^n\Vert_{L^1_{loc}
(\R_+;\,\dot B^{\frac{N}{p_2}-2}_{p_2,1})}
\lesssim
\Vert{\mathcal{Q}}f\Vert_{L^1_{loc}
(\R_+;\,\dot B^{\frac{N}{p_2}-2}_{p_2,1})}.
\end{aligned}
$$
Now, owing  to [\cite{AH}, Theorem 1.1], we deduce that  system
$(\rm\widetilde{MHD})$ with the initial data $(a_0^n,u_0^n,B_0^n,f^n)$
admits a unique local in time solution  $(a^n,u^n,\nabla\Pi^n,B^n)$
verifying
$$
\begin{aligned}
a^n\in C([0,T^n);&\,H^{s+1}(\R^N)),\; u^n,\,B^n\in
C([0,T^n);\,H^s(\R^N))\cap\widetilde L^1_{T^n}(H^{s+2})
\\&
\mbox{and}\quad\nabla\Pi^n\in
L^1([0,T^n);\,H^s(\R^N))\quad\mbox{with}\quad s>\frac{N}{2}-1.
\end{aligned}
$$
\subsubsection*{Estimates of the regularized solution.}
Let $T\in[0,+\infty]$  be defined as
$\displaystyle\inf_{n\in\N}T^n.$ Our first goal is to prove that
$T>0$ such that $(a^n,u^n,\nabla\Pi^n,B^n)$ belongs to and is
uniformly bounded in the space
$$
\begin{aligned}
E_T=
\Big(\widetilde{L}^\infty_T(\dot B^{\frac{N}{p_1}}_{p_1,1}\big)\Big)
&\times\Big(L^1_T\big(\dot B^{\frac{N}{p_2}+1}_{p_2,1}\big)
\cap\widetilde{L}^\infty_T\big(\dot
B^{\frac{N}{p_2}-1}_{p_2,1}\big)\Big)
&\times L^1_T\big(\dot B^{\frac{N}{p_2}-1}_{p_2,1}\big)
\\&
\times\Big(L^1_T\big(\dot B^{\frac{N}{p_2}+1}_{p_2,1}\big)
\cap\widetilde{L}^\infty_T\big(\dot B^{\frac{N}{p_2}-1}_{p_2,1}\big)\Big).
\end{aligned}
$$
Let $(u_L^n,\Pi_L^n)$ be a solution of the following non-stationary 
Stokes system
$$
(\rm{L})\quad\left\{\begin{array}{rl}
&\partial_t u_L^n -\mu^1\Delta u_L^n+\nabla\Pi_L^n=f^n
\medskip\\
&\partial_t B_L^n -\sigma^1\Delta B_L^n=0
\medskip\\
&{\mathop{\rm div}} \,u_L^n ={\mathop{\rm div}} \,B_L^n=0
\medskip\\
&({u_L^n},{B_L^n})_{|t=0}=(u_0^n,B^n_0).
\end{array}
\right.
$$
By construction, $u_0^n,\,B^n_0\in\dot B^{\frac{N}{p_2}-1}_{p_2,1}
\cap H^s$ and
$f^n\in L^1_{loc}(\R_+;\,\dot B^{\frac{N}{p_2}-1}_{p_2,1}\cap H^s).$
So following Proposition 2.3 from \cite{DAN}, we have
$(u^n_L,\nabla\Pi^n_L,B^n_L)\in
L^\infty_t(\dot B^{\frac{N}{p_2}-1}_{p_2,1}\cap H^s)
\times L^1_t(\dot B^{\frac{N}{p_2}-1}_{p_2,1}\cap H^{s})
\times L^\infty_t(\dot B^{\frac{N}{p_2}-1}_{p_2,1}\cap H^s)$ and
moreover $u^n_L,\,B^n_L\in L^1_t(\dot B^{\frac{N}{p_2}+1}_{p_2,1})$
for all $t>0.$\\
Let $u^n=u_L^n+\overline u^n,$
$\nabla\Pi^n=\nabla\Pi^n_L+\nabla\overline\Pi^n$ and
$B^n=B_L^n+\overline B^n.$ Then
$$
\begin{aligned}
(a^n,\overline u^n,\nabla\overline\Pi^n,\overline B^n)
&\in
C\big(([0,T^n);\,H^{s+1}(\R^N)\big)
\times\big(C[0,T^n);\,H^s(\R^N)\big)
\\&
\times L^1_{T^n}\big(H^s(\R^N)\big)
\times C\big([0,T^n);\,H^s(\R^N)\big)
\end{aligned}
$$
and verifies
$$
(\rm{NL})\;\left\{\begin{array}{rl}
&\hspace{-0,5cm}\partial_ta^n+u^n\cdot\nabla a^n=0
\medskip\\
&\hspace{-0,5cm}\partial_t\overline u^n+u^n\cdot\nabla\overline
u^n-\mu^1\Delta\overline u^n+\nabla\overline\Pi^n
=H(a^n,u^n,\nabla\Pi^n,B^n)
\medskip\\
&\hspace{-0,5cm}\partial_t\overline B^n+u^n\cdot\nabla\overline
B^n-\sigma^1\Delta\overline B^n
=-{\mathop{\rm div}}\Big[\big(\widetilde\sigma(a^n)-\sigma^1\big)\nabla B^n\Big]
+B^n\cdot\nabla u^n
\\&
\hspace{5cm}-u^n\cdot\nabla B^n_L
\medskip\\
&\hspace{-0,5cm}{\mathop{\rm div}}\,\overline u^n={\mathop{\rm div}}\,\overline B^n=0
\medskip\\
&\hspace{-0,5cm}(a^n,\overline u^n,\overline B^n)_{|t=0}=(a^n_0,0,0),
\end{array}
\right.
$$
where
$$
\begin{aligned}
H(a^n&,u^n,\nabla\Pi^n,B^n)
=
-u^n\cdot\nabla u^n_L+a^n(\mu^1\Delta u^n-\nabla\Pi^n) 
\\&
+2(1+a^n){\mathop{\rm div}}
\Big\{\big(\widetilde\mu(a^n)-\mu^1\big){\mathcal{M}}^n\Big\}
+(1+a^n)\big(B^n\cdot\nabla B^n-\frac{1}{2}\nabla{B^n}^2\big)
\end{aligned}
$$
with
${\mathcal{M}}^n=\frac{1}{2}(\nabla u^n+^t\nabla u^n).$ We find that
$
(a^n,\overline u^n,\nabla\overline\Pi^n,\overline B^n)
$
belongs to
$
E_{T^n}.
$
by following the arguments as in \cite{HA}.

Now we are in a position to prove that $T>0$ such that
$(a^n,u^n,\nabla\Pi^n,B^n)$
is bounded in $E_{T}.$\\
In what follows, we will use the notation
$$
U^n(t):=
\Vert\overline u^n\Vert_{\widetilde{L}^{\infty}_{t}
(\dot B^{\frac{N}{p_2}-1}_{p_2,1})}
+
\Vert\overline u^n\Vert_{L^1_{t}(\dot B^{\frac{N}{p_2}+1}_{p_2,1})}
+
\Vert\nabla\overline\Pi^n\Vert_{L^1_{t}
(\dot B^{\frac{N}{p_2}-1}_{p_2,1})}
$$
and
$$
B^n(t):=
\Vert\overline B^n\Vert_{\widetilde{L}^{\infty}_{t}(\dot 
B^{\frac{N}{p_2}-1}_{p_2,1})}
+
\Vert\overline B^n\Vert_{L^1_{t}(\dot B^{\frac{N}{p_2}+1}_{p_2,1})}.
$$
Since $\frac{1}{p_1}\leq \frac{1}{N}+{1\over p_2},$
then according to Proposition \ref{eqtransport}, we have
$$
\begin{aligned}
\Vert a^n\Vert_{\widetilde{L}^{\infty}_{T^n}
(\dot B^{\frac{N}{p_1}}_{p_1,1})}
&\leq
e^{C\Vert u^n\Vert_{L^1_{T^n}(\dot B^{\frac{N}{p_2}+1}_{p_2,1})}}
\Vert a_0^n\Vert_{\dot B^{\frac{N}{p_1}}_{p_1,1}}
\\&
\lesssim
e^{C\Vert u^n\Vert_{L^1_{T^n}(\dot B^{\frac{N}{p_2}+1}_{p_2,1})}}
\Vert a_0\Vert_{\dot B^{\frac{N}{p_1}}_{p_1,1}}.
\end{aligned}
$$
Moreover, Proposition \ref{eq.stokes} implies that
$$
U^n(T^n)
\leq
Ce^{C\Vert\nabla u^n\Vert_{L^1_{T^n}(\dot B^{\frac{N}{p_2}}_{p_2,1})}}
\Vert H(a^n,u^n,\nabla\Pi^n,B^n)\Vert_{L^1_{T^n}
(\dot B^{\frac{N}{p_2}-1}_{p_2,1})}.
$$
Since ${1\over p_1}+{1\over p_2}>{1\over N},$ then the inequality 
(\ref{p_1p_2}) implies that
\begin{equation}\label{1es-Existence}
\Vert
a^n(\mu^1\Delta u^n-\nabla\Pi^n)
\Big\Vert_{L^1_{T^n}(\dot B^{\frac{N}{p_2}-1}_{p_2,1})}
\lesssim
\Vert a^n\Vert_{L^\infty_{T^n}(\dot B^{\frac{N}{p_1}}_{p_1,1})}
\Big(\Vert u^n\Vert_{L^1_{T^n}(\dot B^{\frac{N}{p_2}+1}_{p_2,1})}
+\Vert\nabla\Pi^n\Vert_{L^1_{T^n}
(\dot B^{\frac{N}{p_2}-1}_{p_2,1})}\Big).
\end{equation}
From the Bernstein inequality, (\ref{p_1p_2}) and a classical 
interpolation argument, we may infer that
\begin{equation}\label{1es-Exis}
\begin{aligned}
\Vert u^n\cdot\nabla u^n_L\Vert_{L^1_{T^n}
(\dot B^{\frac{N}{p_2}-1}_{p_2,1})}
&\lesssim
\Vert u^n\otimes  u^n_L\Vert_{L^1_{T^n}
(\dot B^{\frac{N}{p_2}}_{p_2,1})}
\\&
\lesssim
\Vert u^n\Vert_{L^2_{T^n}(\dot B^{\frac{N}{p_2}}_{p_2,1})}
\Vert u^n_L\Vert_{L^2_{T^n}(\dot B^{\frac{N}{p_2}}_{p_2,1})}.
\end{aligned}
\end{equation}
Since $p_1\leq p_2$ and ${1\over p_1}+{1\over p_2}>{1\over N},$ then the Bernstein inequality, estimate 
(\ref{p_1p_2}),  and Taylor's formula imply that
$$
\begin{aligned}
\Big\Vert(1+a^n){\mathop{\rm
div}}\big\{\big(\widetilde\mu(a^n)&-\mu^1\big)
{\mathcal{M}}^n\big\}\Big\Vert_{L^1_{T^n}
(\dot B^{\frac{N}{p_2}-1}_{p_2,1})}
\lesssim
\big(1+\Vert a^n\Vert_{L^\infty_{T^n}(\dot B^{\frac{N}{p_1}}_{p_1,1})}\big)
\\&
\times
\big\Vert\big(\widetilde\mu(a^n)-\mu^1\big){\mathcal{M}}^n
\big\Vert_{L^1_{T^n}(\dot B^{\frac{N}{p_2}}_{p_2,1})}
\\&
\lesssim
\big(1+\Vert a^n\Vert_{L^\infty_{T^n}
(\dot B^{\frac{N}{p_1}}_{p_1,1})}\big)
\Vert\widetilde\mu(a^n)-\mu^1\Vert_{L^\infty_{T^n}
(\dot B^{\frac{N}{p_1}}_{p_1,1})}
\Vert u^n\Vert_{L^1_{T^n}(\dot B^{\frac{N}{p_2}+1}_{p_2,1})}
\\&
\lesssim
\big(1+\Vert a^n\Vert_{L^\infty_{T^n}
(\dot B^{\frac{N}{p_1}}_{p_1,1})}\big)
\Vert a^n\Vert_{\widetilde L^\infty_{T^n}
(\dot B^{\frac{N}{p_1}}_{p_1,1})}
\Vert u^n\Vert_{L^1_{T^n}(\dot B^{\frac{N}{p_2}+1}_{p_2,1})}
\end{aligned}
$$
and
$$
\begin{aligned}
\Big\Vert(1+a^n)\big(B^n\cdot\nabla B^n&-\frac{1}{2}\nabla{B^n}^2\big)
\Big\Vert_{L^1_{T^n}(\dot B^{\frac{N}{p_2}-1}_{p_2,1})}
\lesssim
\Big(1+\Vert a^n\Vert_{L^\infty_{T^n}
(\dot B^{\frac{N}{p_1}}_{p_1,1})}\Big)
\\&
\times
\Big(\Vert B^n\otimes B^n\Vert_{L^1_{T^n}
(\dot B^{\frac{N}{p_2}}_{p_2,1})}
+\Vert {B^n}^2\Vert_{L^1_{T^n}(\dot B^{\frac{N}{p_2}}_{p_2,1})}\Big)
\\&
\lesssim
\big(1+\Vert a^n\Vert_{L^\infty_{T^n}
(\dot B^{\frac{N}{p_1}}_{p_1,1})}\big)
\Vert B^n\Vert_{L^\infty_{T^n}(\dot B^{\frac{N}{p_2}-1}_{p_2,1})}
\Vert B^n\Vert_{L^1_{T^n}(\dot B^{\frac{N}{p_2}+1}_{p_2,1})}.
\end{aligned}
$$
For $\overline B^n,$ we have
$$
\partial_t\overline B^n+u^n\cdot\nabla\overline
B^n-\sigma^1\Delta\overline B^n
=-{\mathop{\rm div}}\Big[\big(\widetilde\sigma(a^n)-\sigma^1\big)\nabla 
B^n\Big]
+\overline B^n\cdot\nabla u^n+B^n_L\cdot\nabla u^n-u^n\cdot\nabla 
B^n_L,
$$
By Proposition \ref{eq.stokes} and inequalities
(\ref{p_1p_2}) and (\ref{produit3}) that for $t\in[0,T^n]$
$$
\begin{aligned}
B^n(t)
\leq
Ce^{C\Vert u^n\Vert_{L^1_{t}(\dot B^{\frac{N}{p_2}+1}_{p_2,1})}}
\Big\{\Vert B^n_L\Vert_{L^2_{t}(\dot B^{\frac{N}{p_2}}_{p_2,1})}
&\Vert u^n\Vert_{L^2_{t}(\dot B^{\frac{N}{p_2}}_{p_2,1})}
+\Vert a^n\Vert_{\widetilde L^\infty_{t}
(\dot B^{\frac{N}{p_1}}_{p_1,1})}
\Vert B^n\Vert_{L^1_{t}(\dot B^{\frac{N}{p_2}+1}_{p_2,1})}
\\&
+\Vert u^n\Vert_{L^\infty_{t}(\dot B^{\frac{N}{p_2}-1}_{p_2,1})}
\Vert B^n_L\Vert_{L^1_{t}(\dot B^{\frac{N}{p_2}+1}_{p_2,1})}\Big\}.
\end{aligned}
$$
So, by interpolation, we have
$$
\Vert v\Vert_{L^2_{t}(\dot B^{\frac{N}{p_2}}_{p_2,1})}
\leq
\Vert v\Vert_{L^1_{t}(\dot B^{\frac{N}{p_2}+1}_{p_2,1})}^{1\over 2}
\Vert v\Vert_{L^\infty_{t}(\dot B^{-1+\frac{N}{p_2}}_{p_2,1})}^{1\over 2}
\hspace{1cm}\forall\quad v\in L^1_{t}(\dot B^{\frac{N}{p_2}+1}_{p_2,1})
\cap
L^\infty_{t}(\dot B^{\frac{N}{p_2}-1}_{p_2,1}),
$$
thus
$$
\begin{aligned}
B^n(t)
&\lesssim
e^{C\Vert u^n\Vert_{L^1_{t}(\dot B^{\frac{N}{p_1}+1}_{p_2,1})}}
\Bigg\{
\Vert B^n_L\Vert_{L^\infty_{t}
(\dot B^{\frac{N}{p_2}-1}_{p_2,1})}^{\frac{1}{2}}
\Vert B^n_L\Vert_{L^1_{t}(\dot 
B^{\frac{N}{p_2}+1}_{p_2,1})}^{\frac{1}{2}}
\Vert u^n\Vert_{L^\infty_{t}(\dot 
B^{\frac{N}{p_2}-1}_{p_2,1})}^{\frac{1}{2}}
\Vert u^n\Vert_{L^1_{t}(\dot 
B^{\frac{N}{p_2}+1}_{p_2,1})}^{\frac{1}{2}}
\\&
+\Vert a^n\Vert_{\widetilde L^\infty_{t}
(\dot B^{\frac{N}{p_1}}_{p_1,1})}
\Vert B^n\Vert_{L^1_{t}(\dot B^{\frac{N}{p_2}+1}_{p_2,1})}
+\Vert u^n\Vert_{L^\infty_{t}(\dot B^{\frac{N}{p_2}-1}_{p_2,1})}
\Vert B^n_L\Vert_{L^1_{t}(\dot B^{\frac{N}{p_2}+1}_{p_2,1})}
\Bigg\}.
\end{aligned}
$$
In the same manner, we have
\begin{equation}\label{Vitesse}
\begin{aligned}
U^n(t)
&\lesssim
 e^{C\Vert u^n\Vert_{L^1_{t}(\dot B^{\frac{N}{p_2}+1}_{p_2,1})}}
\Bigg[
\Vert u^n_L\Vert_{L^\infty_{t(\dot B^{\frac{N}{p_2}-1}_{p_2,1})}}
^{\frac{1}{2}}
\Vert u^n_L\Vert_{L^1_{t}(\dot B^{\frac{N}{p_2}+1}_{p_2,1})}
^{\frac{1}{2}}
\Vert u^n\Vert_{L^\infty_{t}(\dot 
B^{\frac{N}{p_2}-1}_{p_2,1})}^{\frac{1}{2}}
\Vert u^n\Vert_{L^1_{t}(\dot 
B^{\frac{N}{p_2}+1}_{p_2,1})}^{\frac{1}{2}}
\\&
+\Vert a^n\Vert_{\widetilde L^\infty_{t}
(\dot B^{\frac{N}{p_1}}_{p_1,1})}
\Big(1+\Vert a^n\Vert_{L^\infty_{t}
(\dot B^{\frac{N}{p_1}}_{p_1,1})}\Big)
\Big(\Vert u^n\Vert_{L^1_{t}(\dot B^{\frac{N}{p_2}+1}_{p_2,1})}
+\Vert\nabla\Pi^n\Vert_{L^1_{t}(\dot B^{\frac{N}{p_2}-1}_{p_2,1})}
\\&\hspace{7cm}
+
\Vert B^n\Vert_{L^\infty_{t}(\dot B^{\frac{N}{p_2}-1}_{p_2,1})}
\Vert B^n\Vert_{L^1_{t}(\dot B^{\frac{N}{p_2}+1}_{p_2,1})}\Big)
\Bigg].
\end{aligned}
\end{equation}
Let $\zeta$ be a small positive real number. Then there exists
$T_1>0$ such that
\begin{equation}\label{2es-Existence}
\Vert(u_L,B_{L})\Vert_{L^1_{T_1}(\dot B^{\frac{N}{p_2}+1}_{p_2,1})}
+\Vert\nabla\Pi_L\Vert_{L^1_{T_1}(\dot B^{\frac{N}{p_2}-1}_{p_2,1})}
\leq \zeta
\end{equation}
and (see Proposition 2.3 of \cite{DAN})
$$
\Vert u_L\Vert_{\widetilde{L}^{\infty}_{T_1}
(\dot B^{\frac{N}{p_2}-1}_{p_2,1})}
\leq
\Vert u_0\Vert_{\dot B^{\frac{N}{p_2}-1}_{p_2,1}}
+\Vert{\mathcal{P}}f\Vert_{L^1_{T_1}(\dot B^{\frac{N}{p_2}-1}_{p_2,1})}
:=U_0.
$$
Consequently we have
\begin{equation}\label{5es-Existence}
\Vert u_L^n\Vert_{L^1_{T_1}(\dot B^{\frac{N}{p_2}+1}_{p_2,1})}
+\Vert\nabla\Pi_L^n\Vert_{L^1_{T_1}(\dot B^{\frac{N}{p_2}-1}_{p_2,1})}
\leq
C\zeta\;\quad\mbox{et}\quad\;
\Vert u_L^n\Vert_{\widetilde{L}^{\infty}_{T_1}
(\dot B^{\frac{N}{p_2}-1}_{p_2,1})}
\leq
CU_0
\end{equation}
and
\begin{equation}\label{6es_Existence}
\Vert B_L^n\Vert_{L^1_{T_1}(\dot B^{\frac{N}{p_2}+1}_{p_2,1})} 
\leq
C\zeta\;\quad\mbox{and}\quad\; \Vert
B_L^n\Vert_{\widetilde{L}^{\infty}_{T_1}
(\dot B^{\frac{N}{p_2}-1}_{p_2,1})} 
\leq 
C\Vert B_0\Vert_{\dot B^{\frac{N}{p_2}-1}_{p_2,1}}.
\end{equation}
In the following we can suppose that $T^n\leq T_1$, otherwise we take a
smaller $T^n.$ Let $t\leq T^n,$ then
$$
\begin{aligned}
B^n(t)
&\leq
Ce^{C\Big(\zeta+\Vert\overline u^n\Vert_{L^1_{t}
(\dot B^{\frac{N}{p_2}+1}_{p_2,1})}\Big)}
\bigg\{
{\zeta}^{\frac{1}{2}}
\Big(U_0+\Vert\overline u^n\Vert_{L^\infty_{t}
(\dot B^{\frac{N}{p_2}-1}_{p_2,1})}\Big)^{\frac{1}{2}}
\Big(\zeta+\Vert\overline u^n\Vert_{L^1_{t}
(\dot B^{\frac{N}{p_2}+1}_{p_2,1})}\Big)^{\frac{1}{2}}
\\&
+\Vert a^n\Vert_{\widetilde L^\infty_{t}
(\dot B^{\frac{N}{p_1}}_{p_1,1})}
\Big(\zeta+\Vert\overline B^n\Vert_{L^1_{t}
(\dot B^{\frac{N}{p_2}+1}_{p_2,1})}\Big)
+\zeta
\Big(U_0+\Vert\overline u^n\Vert_{L^\infty_{t}
(\dot B^{\frac{N}{p_2}-1}_{p_2,1})}\Big)\bigg\}
\end{aligned}
$$
and
\begin{equation}\label{7es_Existence}
\Vert a^n\Vert_{\widetilde{L}^{\infty}_{t}
(\dot B^{\frac{N}{p_1}}_{p_1,1})}
\leq
Ce^{{C}{}
\Big(\zeta+\Vert\overline u^n\Vert_{L^1_{t}
(\dot B^{\frac{N}{p_2}+1}_{p_2,1})}\Big)}
\Vert a_0\Vert_{\dot B^{\frac{N}{p_1}}_{p_1,1}}.
\end{equation}
Let $T_2\leq T^n$ such that
\begin{equation}\label{temps d'existence}
\exp\Big({C}{}\big(\zeta+\Vert\overline u^n\Vert_{L^1_{T_2}
(\dot B^{\frac{N}{p_2}+1}_{p_2,1})}\big)\Big)< 2.
\end{equation}
So if
$$
16C^2\Vert  a_0\Vert_{\dot B^{\frac{N}{p_1}}_{p_1,1}}
\leq
1,
$$
then
\begin{equation}\label{transport}
\Vert a^n\Vert_{\widetilde{L}^{\infty}_{T_1}
(\dot B^{\frac{N}{p_1}}_{p_1,1})}
\leq
2C\Vert a_0\Vert_{\dot B^{\frac{N}{p_1}}_{p_1,1}}
\end{equation}
and
\begin{equation}\label{poincare}
\begin{aligned}
{\mathcal{B}}^n(T_2)
\leq
4C\bigg\{{\zeta}^{\frac{1}{2}}
\Vert B_0\Vert_{\dot B^{\frac{N}{p_2}-1}_{p_2,1}}^{\frac{1}{2}}
\Big(U_0&+\Vert\overline u^n\Vert_{L^\infty_{t}
(\dot B^{\frac{N}{p_2}-1}_{p_2,1})}\Big)^{\frac{1}{2}}
\Big(\zeta+\Vert\overline u^n\Vert_{L^1_{t}
(\dot B^{\frac{N}{p_2}+1}_{p_2,1})}\Big)^{\frac{1}{2}}
\\&
+2C\zeta\Vert a_0\Vert_{\dot B^{\frac{N}{p}}_{p,1}}
+\zeta\Big(U_0+\Vert\overline u^n\Vert_{L^\infty_{t}
(\dot B^{\frac{N}{p_2}-1}_{p_2,1})}\Big)\bigg\}.
\end{aligned}
\end{equation}
Using inequalities (\ref{Vitesse}) and
(\ref{poincare}) satisfied by $B^n=B_{L}^n+\overline{B}^n$, we obtain 
that
$$
\begin{aligned}
U^n(T_2)
\leq
C\bigg\{
\zeta\big(U^n(T_2)&+U_0\big)
+2C\Vert a_0\Vert_{\dot B^{\frac{N}{p_1}}_{p_1,1}}
\Big(1+2C\Vert a_0\Vert_{\dot B^{\frac{N}{p_1}}_{p_1,1}}\Big)
\Big(\zeta+U^n(T_2)\Big)
\\&
+\zeta\Vert B_0\Vert_{\dot B^{\frac{N}{p_2}-1}_{p_2,1}}
\Big(1+2C\Vert a_0\Vert_{\dot B^{\frac{N}{p_1}}_{p_1,1}}\Big)
\Big(U_0^2+\zeta^2+U^n(T_2)^2\Big)\bigg\}.
\end{aligned}
$$
Using (\ref{transport}) and the smallness of $a_{0},$ we obtain for
$\zeta$ small enough,
\begin{equation}\label{8es_Existence}
U^n(T_2)
\leq
\zeta C\Big(U_0,\Vert a_0\Vert_{\dot B^{\frac{N}{p_1}}_{p_1,1}},
\Vert B_0\Vert_{\dot B^{\frac{N}{p_2}-1}_{p_2,1}}\Big).
\end{equation}
Taking $\zeta$ small enough we observe that inequality
(\ref{temps d'existence}) is satisfied. Consequently,  a standard
argument then yields that $T_2=T^n.$ The same type of
reasoning allows one to show that $T^{n}=T^{1},$ with uniform 
control.\\
We give in what follows a precise estimate of the pressure term.
Namely, we prove the following
\begin{lem}\label{malanelabou}{\it
Let $0<\eta<\inf(1,\frac{2N}{p_2})$ be such that 
$\frac{1}{N}+\frac{\eta}{N}<\frac{1}{p_1}+\frac{1}{p_2}.$
Then
$\nabla\Pi^n$ is uniformly bounded in
$L^{\frac{2}{2-\eta}}_{T_1}(\dot B^{\frac{N}{p_2}-1-\eta}_{p_2,1}).$}
\end{lem}
\noindent {\bf Proof.} Applying the divergence operator to the
equation containing the pressure term, we obtain
$$
\begin{aligned}
{\mathop{\rm div}}\Big(\big(1+a^n\big)\nabla\Pi^n\Big)=
{\mathop{\rm div}}\bigg\{&\Big(1+a^n\Big)\Big({\mathop{\rm div}}
\big\{\widetilde\mu(a^n){\mathcal{M}}^n\big\}+B^n\cdot\nabla B^n
-\frac{1}{2}\nabla{B^n}^2\Big)
\\&
+{\mathcal{Q}}f^n-u^n\cdot\nabla u^n\bigg\}.
\end{aligned}
$$
By construction of $f^n$ and by interpolation, we have
that ${\mathcal{Q}}f^n$ is uniformly bounded in
$L^{\frac{2}{2-\eta}}_{T_1}(\dot B^{\frac{N}{p_2}-1-\eta}_{p_2,1}).$
By interpolation, we have that  $u^n$ is uniformly bounded
in $L^{\frac{2}{1-\eta}}_{T_1}(\dot B^{\frac{N}{p_2}-\eta}_{p_2,1}).$ 
Since $\eta<\frac{2N}{p_2},$
 inequality (\ref{produit3}) implies the estimate
$$
\begin{aligned}
\Vert u^n\cdot\nabla u^n\Vert_{L^{\frac{2}{2-\eta}}_{T_1}
(\dot B^{\frac{N}{p_2}-1-\eta}_{p_2,1})}
&\lesssim
\Vert u^n\otimes u^n\Vert_{L^{\frac{2}{2-\eta}}_{T_1}
(\dot B^{\frac{N}{p_2}-\eta}_{p_2,1})}
\\&
\lesssim
\Vert u^n\Vert_{L^{\frac{2}{1-\eta}}_{T_1}
(\dot B^{\frac{N}{p_2}-\eta}_{p_2,1})}
\Vert u^n\Vert_{L^2_{T_1}(\dot B^{\frac{N}{p_2}}_{p_2,1})},
\end{aligned}
$$
which shows that $u^n\cdot\nabla u^n$ is
uniformly bounded in $L^{\frac{2}{2-\eta}}_{T_1}
(\dot B^{\frac{N}{p_2}-1-\eta}_{p_2,1}).$
In the same way $\big(1+a^n\big){\mathop{\rm div}}
\big\{\widetilde\mu(a^n){\mathcal{M}}^n\big\},$ for $p_1\leq p_2,$
and $\frac{1}{p_1}+\frac{1}{p_2}>\frac{1}{N},$
the Bernstein inequality  and (\ref{p_1p_2}) imply that
${\mathop{\rm div}}\big\{\widetilde\mu(a^n){\mathcal{M}}^n\big\}$ is 
uniformly bounded
in $L^1_{T_1}(\dot B^{\frac{N}{p_2}-1}_{p_2,1}) \cap L^2_{T_1}
(\dot B^{\frac{N}{p_2}-2}_{p_2,1}).$ 
So, by an interpolation argument, we
obtain that ${\mathop{\rm div}}\big\{\widetilde\mu(a^n)
{\mathcal{M}}^n\big\}$ is uniformly
bounded in $L^{\frac{2}{2-\eta}}_{T_1}
(\dot B^{\frac{N}{p_2}-1-\eta}_{p_2,1}).$ Since 
$\frac{1}{N}+\frac{\eta}{N}<\frac{1}{p_1}+\frac{1}{p_2},$
the inequality (\ref{p_1p_2}) implies that\\
$\big(1+a^n\big){\mathop{\rm div}}\big\{\widetilde\mu(a^n) 
{\mathcal{M}}^n\big\}$ is uniformly
bounded in  $L^{\frac{2}{2-\eta}}_{T_1}
(\dot B^{\frac{N}{p_2}-1-\eta}_{p_2,1})$ in the same way as for
$\big(1+a^n\big)\big(B^n\cdot\nabla B^n-\frac{1}{2}\nabla B^n\big).$
So $\nabla\Pi^n$ is also uniformly bounded because we have
$\Vert a^n\Vert_{\widetilde{L}^{\infty}_{T_1}
(\dot B^{\frac{N}{p_1}}_{p_1,1})}
\leq
2C\Vert a_0\Vert_{\dot B^{\frac{N}{p_1}}_{p_1,1}}<<1.$
$\square$

By the construction of the time of existence, then 
$T_1=\infty,$  provided that
$$
\Vert u_0\Vert_{\dot B^{\frac{N}{p_2}-1}_{p_2,1}}
+\Vert B_0\Vert_{\dot B^{\frac{N}{p_2}-1}_{p_2,1}}
+\Vert f\Vert_{L^1(\R_+;\,\dot B^{\frac{N}{p_2}-1}_{p_2,1})}
\leq
c'\inf(\mu^1,\sigma^1).
$$
\subsubsection*{Passage to the limit.}
Let us note first that by construction of
$(u^n_0,f^n),$ the sequence $(u^n_L,\nabla\Pi^n_L,B^n_L)$ converges
strongly to the solution $(u_L,\nabla\Pi_L,B_L)$ of the  system
$(L).$ However, to show that the weak limit of $(a^n,\overline
u^n,\nabla\overline\Pi^n,\overline B^n)$ is a solution to the
system $(NL),$ we need to use some compactness arguments.\\
We have already established that $(a^n,\overline
u^n,\nabla\overline\Pi^n,\overline B^n)$ is uniformly  bounded in
$$
\widetilde{L}^\infty_{T_1}(\dot B^{\frac{N}{p_1}}_{p_1,1})
\times\widetilde{L}^\infty_{T_1}(\dot B^{\frac{N}{p_2}-1}_{p_2,1})
\cap L^1_{T_1}(\dot B^{\frac{N}{p_2}+1}_{p_2,1})
\times L^1_{T_1}(\dot B^{\frac{N}{p_2}-1}_{p_2,1})
\times\widetilde{L}^\infty_{T_1}(\dot B^{\frac{N}{p_2}-1}_{p_2,1})
\cap L^1_{T_1}(\dot B^{\frac{N}{p_2}+1}_{p_2,1}),
$$
Moreover $\nabla\Pi^n$ is uniformly bounded in
$L^{\frac{2}{2-\eta}}_{T_1}(\dot B^{\frac{N}{p_2}-1-\eta}_{p_2,1}).$\\
So, in order to use the Ascoli theorem, it suffices to  estimate the
time derivative of $a^n,$ $\overline u^n$ and $\overline B^n$ 
(see for example \cite{DANC}).
Following the proof of Lemma \ref{malanelabou}, the following lemma is shown to hold true.
\begin{lem}\label{Th-Ascoli}
{\it\quad\\
(i) The sequence $(\partial_ta^n)_{n\in\N}$ is uniformly bounded in
$L^2_{T_1}(\dot B^{\frac{N}{p_1}-1}_{p_1,1}).$\\
(ii) The sequence $(\partial_t\overline u^n)_{n\in\N}$ is uniformly
bounded in $L^{\frac{2}{2-\eta}}_{T_1}(\dot 
B^{\frac{N}{p_2}-1-\eta}_{p_2,1})$
for \\
$0<\eta<\inf(1,\frac{2N}{p_2})$ and 
$\frac{1}{N}+\frac{\eta}{N}<\frac{1}{p_1}+\frac{1}{p_2}.$\\
(iii) The sequence $(\partial_t\overline B^n)_{n\in\N}$ is uniformly
bounded in $L^{\frac{2}{2-\eta}}_{T_1}(\dot 
B^{\frac{N}{p_2}-1-\eta}_{p_2,1})$
for \\
$0<\eta<\inf(1,\frac{2N}{p_2})$ and 
$\frac{1}{N}+\frac{\eta}{N}<\frac{1}{p_1}+\frac{1}{p_2}.$}
\end{lem}

From the above lemma, the Cauchy-Schwarz inequality and
H\"older's inequality, we deduce the following corollary.
\begin{corol}\label{Lemma-2}\quad\\
(i) The sequence $(a^n)_{n\in\N}$ is uniformly bounded in
$C^{\frac{1}{2}}\Big([0,T_1];\,\dot 
B^{\frac{N}{p_1}-1}_{p_1,1}\Big).$\\
(ii) The sequence $(\overline u^n)_{n\in\N}$ is uniformly bounded in
$C^{\frac{\eta}{2}}\Big([0,T_1];\,\dot 
B^{\frac{N}{p_2}-1-\eta}_{p_2,1}\Big)$
for all $\eta$ belonging to $]0,\inf(1,\frac{2N}{p_2})[$ and
$\frac{1}{N}+\frac{\eta}{N}<\frac{1}{p_1}+\frac{1}{p_2}.$\\
(iii) The sequence $(\overline B^n)_{n\in\N}$ is uniformly bounded
in $C^{\frac{\eta}{2}}\Big([0,T_1];\,\dot 
B^{\frac{N}{p_2}-2}_{p_2,1}\Big)$ for all $\eta$ belonging to
$]0,\inf(1,\frac{2N}{p_2})[$ and
$\frac{1}{N}+\frac{\eta}{N}<\frac{1}{p_1}+\frac{1}{p_2}.$
\end{corol}
We recall that the injection of $\dot B^{s+\varepsilon}_{p\,q,loc}$
in $B^{s}_{p\,q,loc}$ (the inhomogeneous Besov space $B^{s}_{p\,q,loc}$)
is compact for all $\varepsilon>0$ (see for example \cite{RS}).\\
So, there exists a subsequence (which is still denoted by  
$(a^n,\overline u^n,\nabla\overline\Pi^n,\overline B^n)$)  which converges to
$(a,\overline u,\nabla\overline\Pi,\overline B).$ Consequently,
$(a,u,\nabla\Pi,B)$ is a solution of the $(\rm\widetilde{MHD})$ system
belonging to
$$
\widetilde{L}^\infty_{T_1}(\dot B^{\frac{N}{p_1}}_{p_1,1})
\times\widetilde{L}^\infty_{T_1}(\dot B^{\frac{N}{p_2}-1}_{p_2,1})
\cap L^1_{T_1}(\dot B^{\frac{N}{p_2}+1}_{p_2,1})
\times L^1_{T_1}(\dot B^{\frac{N}{p_2}-1}_{p_2,1})
\times
\widetilde{L}^\infty_{T_1}(\dot B^{\frac{N}{p_2}-1}_{p_2,1})
\cap L^1_{T_1}(\dot B^{\frac{N}{p_2}+1}_{p_2,1}).
$$
Concerning the continuity of $u,$ we have used the fact that
$$
(\rm{H})\;\left\{\begin{array}{rl}
&\hspace{-0,5cm}\partial_tu-\mu^1\Delta u=
H(a,u,\nabla\Pi,B)
\\
&\hspace{-0,5cm}u_{|t=0}=u_0,
\end{array}
\right.
$$
where
$$
\begin{aligned}
H(a,u,\nabla\Pi,B)
=
f-u\cdot\nabla u
&-(1+a)\big(\nabla\Pi+\frac{1}{2}\nabla{B}^2-B\cdot\nabla B\big)
\\&
+2(1+a){\mathop{\rm div}}
\Big\{\big(\widetilde\mu(a)-\mu^1\big){\mathcal{M}}\Big\}
+\mu^1a\Delta u.
\end{aligned}
$$
Since
$$
\begin{aligned}
(a,u,\nabla\Pi,B)\in
\widetilde{L}^\infty_{T_1}(\dot B^{\frac{N}{p_1}}_{p_1,1})
&\times\widetilde{L}^\infty_{T_1}(\dot B^{\frac{N}{p_2}-1}_{p_2,1})
\cap L^1_{T_1}(\dot B^{\frac{N}{p_2}+1}_{p_2,1})
\times L^1_{T_1}(\dot B^{\frac{N}{p_2}-1}_{p_2,1})
\\&
\times
\widetilde{L}^\infty_{T_1}(\dot B^{\frac{N}{p_2}-1}_{p_2,1})
\cap L^1_{T_1}(\dot B^{\frac{N}{p_2}+1}_{p_2,1}),
\end{aligned}
$$
then Proposition \ref{Loi-Pro-Par}, implies that
$H(a,u,\nabla\Pi,B)\in L^1_{T_1}(\dot B^{-1+{N\over p_2}}_{p_2,1}).$
And consequently, Proposition 2.1 \cite{CH}, ensured the continuity
in time of $u,$ in the same way for $B.$
To prove that $a$ is continuous and that the $L^\infty$-norm  is conserved, we use that
$a=a_0\circ\Psi^{-1}$ where $\Psi$ is the flow of $u.$ This
completes the proof of  Theorem \ref{Res-Pri}. 
$\square$

\end{document}